\documentclass[leqno,12pt]{article}  
\usepackage{amsmath}
\usepackage{amssymb}

\usepackage[german,english]{babel}
\usepackage{amsthm}
\usepackage{xypic}  

\author{\textsc{Elmar Grosse-Kl\"onne}}
\date{}

\theoremstyle{plain} 
\newtheorem{satz}{Theorem}[section]  
\newtheorem{kor}[satz]{Corollary}  
\newtheorem{lem}[satz]{Lemma}  
\newtheorem{pro}[satz]{Proposition}  
 
\theoremstyle{remark}

\theoremstyle{definition}

\newcommand{\0}{\ensuremath{\overrightarrow{0}}}

\begin{document}

%

\begin{center}{\bf From pro-$p$ Iwahori-Hecke modules to
    $(\varphi,\Gamma)$-modules, II}\\by Elmar Grosse-Kl\"onne
\end{center}

\begin{abstract} Let ${\mathfrak o}$ be the ring of integers in a finite
  extension field of
${\mathbb Q}_p$, let $k$ be its residue field. Let $G$ be a split reductive group over ${\mathbb
  Q}_p$, let ${\mathcal
  H}(G,I_0)$ be its pro-$p$-Iwahori Hecke ${\mathfrak o}$-algebra. In \cite{dfun} we introduced a general principle how to
  assign to a certain additionally chosen datum $(C^{(\bullet)},\phi,\tau)$ an
  exact functor $M\mapsto{\bf D}(\Theta_*{\mathcal V}_M)$  from finite length ${\mathcal
  H}(G,I_0)$-modules to $(\varphi^r,\Gamma)$-modules. In the present paper we
concretely work out such data $(C^{(\bullet)},\phi,\tau)$ for the classical
matrix groups. We show that the corresponding functor identifies the set of
(standard) supersingular ${\mathcal
  H}(G,I_0)\otimes_{{\mathfrak o}}k$-modules with the set of
$(\varphi^r,\Gamma)$-modules satisfying a certain symmetry condition.
\end{abstract}

\tableofcontents

\section{Introduction}

Let ${\mathfrak o}$ be the ring of integers in a finite
  extension field of
${\mathbb Q}_p$, let $k$ be its residue field. Let $G$ be a split reductive group over ${\mathbb
  Q}_p$, let $T$ be a maximal split torus in $G$, let $I_0$ be a
pro-$p$-Iwahori subgroup fixing a chamber $C$ in the $T$-stable apartment of the semi
simple Bruhat Tits building of $G$. Let ${\mathcal
  H}(G,I_0)$ be the pro-$p$-Iwahori Hecke ${\mathfrak o}$-algebra. Let ${\rm
  Mod}^{\rm fin}({\mathcal
  H}(G,I_0))$ denote the category of ${\mathcal
  H}(G,I_0)$-modules of finite ${\mathfrak o}$-length. From a certain
additional datum $(C^{(\bullet)},\phi,\tau)$ we constructed in \cite{dfun} an exact functor $M\mapsto{\bf D}(\Theta_*{\mathcal V}_M)$ from ${\rm
  Mod}^{\rm fin}({\mathcal
  H}(G,I_0))$ to the category of \'{e}tale $(\varphi^r,\Gamma)$-modules (with
$r\in{\mathbb N}$ depending on $\phi$). For $G={\rm GL}_2({\mathbb Q}_p)$, when
precomposed with the functor of taking $I_0$-invariants, this yields the
functor from smooth ${\mathfrak o}$-torsion representations of ${\rm
  GL}_2({\mathbb Q}_p)$ (or at least from those generated by their $I_0$-invariants) to \'{e}tale $(\varphi,\Gamma)$-modules which plays
a crucial role in Colmez' construction of a $p$-adic local Langlands
correspondence for ${\rm
  GL}_2({\mathbb Q}_p)$. In \cite{dfun} we studied in detail the functor $M\mapsto{\bf
  D}(\Theta_*{\mathcal V}_M)$ when $G={\rm GL}_{d+1}({\mathbb Q}_p)$ for
$d\ge1$. In \cite{koz} the situation has been analysed for $G={\rm SL}_{d+1}({\mathbb Q}_p)$. The purpose of the present paper is to explain how the general
construction of \cite{dfun} can be installed concretely for other classical matrix groups $G$ (as well as for $G$'s of type $E_6$, $E_7$).

Recall that $C^{(\bullet)}=(C=C^{(0)}, C^{(1)},C^{(2)},\ldots)$ is a
minimal gallery, starting at $C$, in the $T$-stable apartment, that $\phi\in
N(T)$ is a 'period' of
$C^{(\bullet)}$ and that $\tau$ is a homomorphism from ${\mathbb
  Z}_p^{\times}$ to $T$, compatible with $\phi$ in a suitable sense. The above
$r\in{\mathbb N}$ is just the length of $\phi$. It turns out that $\tau$ must be a minuscule fundamental
coweight (at least if the underlying root system is simple). Conversely, any minuscule fundamental
coweight $\tau$ can be included into a datum $(C^{(\bullet)},\phi,\tau)$, in
such a way that some power of $\tau$ is a power of $\phi$.

For $G={\rm GL}_{d+1}({\mathbb Q}_p)$ we gave explicit choices of
$(C^{(\bullet)},\phi,\tau)$ with $r=1$ in \cite{dfun} (there are essentially just two choices, and these are 'dual' to each other). In the present paper we work out 'privileged' choices
$(C^{(\bullet)},\phi,\tau)$ for the classical matrix groups, as well as for $G$'s of type $E_6$, $E_7$. We mostly consider
$G$ with connected center $Z$. Our choices of $(C^{(\bullet)},\phi,\tau)$ are
such that $\phi\in N(T)$ projects (modulo $ZT_0$, where $T_0$ denotes the
maximal bounded subgroup of $T$) to the affine Weyl group (viewed as a
subgroup of $N(T)/ZT_0$). In particular, up to modifications by elements of
$Z$ these $\phi$ can also be included into data $(C^{(\bullet)},\phi,\tau)$
for the other $G$'s with the same underlying root system, not necessarily with
connected center. We indicate these modifications along the way. Notice that the $\phi\in N(T)$ considered in \cite{dfun} for $G={\rm
  GL}_{d+1}({\mathbb Q}_p)$ does {\it not} project to the affine Weyl group,
only its $(d+1)$-st power (which is considered here) does so. But since the discussion is essentially
the same, our treatment of the case $A$ here is very brief. 

In either case we work out the behaviour of the functor $M\mapsto{\bf
  D}(\Theta_*{\mathcal V}_M)$ on those ${\mathcal
  H}(G,I_0)_k={\mathcal
  H}(G,I_0)\otimes_{\mathfrak o}k$-modules which we call 'standard
supersingular'. Roughly speaking, these are induced from characters of the
pro-$p$-Iwahori Hecke algebra of the corresponding simply connected group. Each irreducible supersingular ${\mathcal
  H}(G,I_0)_k$-module is contained in (and in 'most' cases is equal to) a standard
supersingular ${\mathcal
  H}(G,I_0)_k$-module (and the very few standard
supersingular ${\mathcal
  H}(G,I_0)_k$-modules which are not irreducible supersingular are easily
identified). We show that our functor induces a bijection between the set of
isomorphism classes (resp. certain packets of such if $G={\rm SO}_{2d+1}({\mathbb Q}_p)$) of standard
supersingular ${\mathcal
  H}(G,I_0)_k$-modules and the set of (isomorphism classes of)
\'{e}tale $(\varphi^r,\Gamma)$-modules over $k_{\mathcal E}=k((t))$ which satisfy a certain symmetry condition (depending on the root system
underlying $G$). They are direct sums of one dimensional \'{e}tale $(\varphi^r,\Gamma)$-modules, their dimension is the $k$-dimension of the
corresponding standard supersingular ${\mathcal H}(G,I_0)_k$-module. 

The interest in \'{e}tale $(\varphi^r,\Gamma)$-modules lies in
their relation with ${\rm Gal}_{{\mathbb Q}_p}$-representations. For any
$r\in{\mathbb N}$ there is an
exact
functor from the category of \'{e}tale $(\varphi^r,\Gamma)$-modules to the category of \'{e}tale
$(\varphi,\Gamma)$-modules (it multiplies the rank by the factor $r$), and by
means of
Fontaine's functor, the
latter one is equivalent with the category of ${\rm Gal}_{{\mathbb
    Q}_p}$-representations.

In \cite{dfun} we also explained that a datum $(C^{(\bullet)},\phi)$ alone,
i.e. without a $\tau$ as above, can be used to define an exact functor $M\mapsto{\bf D}(\Theta_*{\mathcal V}_M)$ from ${\rm
  Mod}^{\rm fin}({\mathcal
  H}(G,I_0))$ to the category of \'{e}tale $(\varphi^r,\Gamma_0)$-modules, where $\Gamma_0$ denotes the maximal pro-$p$-subgroup of $\Gamma\cong{\mathbb
  Z}_p^{\times}$. Such data $(C^{(\bullet)},\phi)$ are not tied to co
minuscule coweights and exist in abundance. We do not discuss them here.

We hope that, beyond its immediate purposes as described above, the present paper may also be a useful reference for explicit descriptions of the pro-$p$-Iwahori algebra (in particular with respect to the various Weyl groups involved) of classical matrix groups other than ${\rm GL}_d$ (we could not find such descriptions in the literature).

The outline is as follows. In section 2 we explain the functor from \'{e}tale $(\varphi^r,\Gamma)$-modules to \'{e}tale
$(\varphi,\Gamma)$-modules, and we introduce the 'symmetric' \'{e}tale
$(\varphi^r,\Gamma)$-modules mentioned above, for each of the root systems $C$,
$B$, $D$ and $A$. In section 3, Lemma \ref{concept}, we discuss the relation between the data
$(C^{(\bullet)},\phi,\tau)$ and minuscule fundamental coweights. Our
discussions of classical matrix groups $G$ in section 4 are just
concrete incarnations of Lemma \ref{concept}, although in neither of these
cases there is a need to make formal
reference to Lemma \ref{concept}. On the other hand, in our discussion
of the cases $E_6$ and $E_7$ in section 5 we do invoke Lemma \ref{concept}. We
tried to synchronize our discussions of the various matrix groups. As a
consequence, arguments repeat themselves, and we do not write them out again and again. In the appendix we record
calculations relevant for the cases $E_6$ and $E_7$, carried out with the help
of the computer algebra system {\it sage}. 

{\it Acknowledgments:} I would like to thank Laurent Berger for a helpful discussion
related to this work. I am very grateful to the referees for the careful reading of the manuscript and the detailed suggestions for improvements.

\section{$(\varphi^r,\Gamma)$-modules}

We often regard elements of ${\mathbb F}_p^{\times}$ as elements of ${\mathbb
  Z}_p^{\times}$ by means of the Teichm\"uller lifting. In ${\rm
  GL}_2({\mathbb Z}_p)$ we define the subgroups$$\Gamma=\left( \begin{array}{cc}
    {\mathbb Z}_p^{\times} & 0 \\
    0 & 1\end{array}\right),\quad\quad\quad\Gamma_{{{0}}}=\left( \begin{array}{cc}
    1+p{\mathbb Z}_p & 0 \\
    0 & 1\end{array}\right),\quad\quad\quad{\mathfrak N}_0=\left( \begin{array}{cc} 1 & {\mathbb Z}_p \\
    0 & 1\end{array}\right)$$and the elements$$\varphi=\left(\begin{array}{cc}p&0\\0&1\end{array}\right),\quad\quad\nu=\left(\begin{array}{cc}1&1\\0&1\end{array}\right),\quad\quad
h(x)=\left(\begin{array}{cc}x&0\\0&x^{-1}\end{array}\right),\quad\quad\gamma(x)=\left( \begin{array}{cc} x & 0 \\ 0 &
    1\end{array}\right)$$
where $x\in
{\mathbb Z}_p^{\times}$. Let ${\mathcal O}_{\mathcal E}^+={\mathfrak o}[[{\mathfrak N}_0]]$ denote the completed group ring of
${\mathfrak N}_0$ over ${\mathfrak o}$. Let ${\mathcal O}_{\mathcal E}$ denote the $p$-adic completion of the
localization of ${{{\mathcal O}_{\mathcal E}^+}}$ with respect to the
complement of $\pi_K{{{\mathcal O}_{\mathcal E}^+}}$, where
$\pi_K\in{\mathfrak o}$ is a uniformizer. In the completed group ring
$k_{\mathcal E}^+=k[[{\mathfrak N}_0]]$ we put $t=[\nu]-1$. Let $k_{\mathcal
  E}={\rm Frac}(k_{\mathcal E}^+)={{{\mathcal O}_{\mathcal
      E}}}\otimes_{{\mathfrak o}}k$. We identify $k_{\mathcal E}^+=k[[t]]$
and $k_{\mathcal E}=k((t))$. For definitions and notational conventions concerning \'{e}tale $\varphi^r$-modules and \'{e}tale
$(\varphi^r,\Gamma)$-modules we refer to \cite{dfun}.
  
Let $r\in{\mathbb N}$. Let ${\bf D}=({\bf D},{\varphi}^r_{{\bf D}})$ be an \'{e}tale $\varphi^r$-module over ${\mathcal
  O}_{\mathcal E}$. For $0\le i\le r-1$ let ${\bf D}^{(i)}={\bf D}$ be a copy of ${\bf D}$. For
$1\le i\le r-1$ define
$\varphi_{\widetilde{\bf D}}:{\bf D}^{(i)}\to{\bf D}^{(i-1)}$ to be the
identity map on ${\bf D}$, and define $\varphi_{\widetilde{\bf D}}:{\bf
  D}^{(0)}\to{\bf D}^{(r-1)}$ to be the structure map ${\varphi}^r_{{\bf D}}$
on ${\bf D}$. Together we obtain a ${\mathbb Z}_p$-linear endomorphism
$\varphi_{\widetilde{\bf D}}$ on $$\widetilde{\bf D}=\bigoplus_{i=0}^{r-1}{\bf
  D}^{(i)}.$$Define an ${\mathcal
  O}_{\mathcal E}$-action on $\widetilde{\bf D}$ by the
formula $x\cdot((d_i)_{0\le i\le r-1})=(\varphi^i_{{\mathcal
  O}_{\mathcal E}}(x)d_i)_{0\le i\le r-1}$. Then the endomorphism $\varphi_{\widetilde{\bf D}}$ of
  $\widetilde{\bf D}$ is semilinear with respect to this ${\mathcal
  O}_{\mathcal E}$-action, hence it defines on $\widetilde{\bf
  D}$ the structure of an \'{e}tale $\varphi$-module over ${\mathcal
  O}_{\mathcal E}$, see \cite{dfun}, section 6.3.

Let $\Gamma'$ be an open subgroup of $\Gamma$, let ${\bf D}$ be an \'{e}tale $(\varphi^r,\Gamma')$-module over ${\mathcal
  O}_{\mathcal E}$. Define an action of $\Gamma'$ on $\widetilde{\bf D}$ by$$\gamma\cdot((d_i)_{0\le i\le r-1})=(\gamma\cdot d_i)_{0\le i\le r-1}.$$

\begin{lem}\label{gamsemilin} The $\Gamma'$-action on $\widetilde{\bf D}$
  commutes with $\varphi_{\widetilde{\bf D}}$ and is semilinear with respect to the ${\mathcal
  O}_{\mathcal E}$-action, hence we obtain on $\widetilde{\bf
  D}$ the structure of an \'{e}tale $(\varphi,\Gamma')$-module over ${\mathcal
  O}_{\mathcal E}$. We thus obtain an exact functor from the category of \'{e}tale $(\varphi^r,\Gamma')$-modules to the category of \'{e}tale $(\varphi,\Gamma')$-modules over ${\mathcal
  O}_{\mathcal E}$.
\end{lem}

{\sc Proof:} This is immediate from the respective properties of the
$\Gamma'$-action on ${\bf
  D}$. \hfill$\Box$\\

The conjugation action of $\Gamma$ on ${\mathfrak N}_0$ allows us to define for each $x\in {\mathbb Z}_p^{\times}$ an automorphism $f\mapsto \gamma(x)f\gamma(x^{-1})$ of $k[[t]]=k_{\mathcal E}^+$.

\begin{lem} For $x\in {\mathbb F}_p^{\times}$ and $m, n\in{\mathbb Z}_{\ge0}$ we have\begin{gather}\gamma(x)t^{np^{rm}}\gamma(x^{-1})-(xt)^{np^{rm}}\in t^{(n+1)p^{rm}}k[[t]].\label{ennoabschied}\end{gather}
\end{lem}

{\sc Proof:} Denote again by $x$ be the representative of $x$ in $[1,p-1]$. We compute $[\nu]^x-1=\sum_{j=1}^x{x\choose j}([\nu]-1)^x=\sum_{j=1}^x{x\choose j}t^j$, hence \begin{align}\gamma(x)t\gamma(x^{-1})-xt=([\nu]^x-1)-xt&\in t^2k[[t]],\notag\\\gamma(x)t^n\gamma(x^{-1})-(xt)^n&\in t^{n+1}k[[t]],\notag\\\gamma(x)t^{np^{rm}}\gamma(x^{-1})-(xt)^{np^{rm}}=(\gamma(x)t^{n}\gamma(x^{-1})-(xt)^{n})^{p^{rm}}&\in t^{(n+1)p^{rm}}k[[t]].\notag\end{align}
 \hfill$\Box$\\

\begin{lem}\label{flippa} (a) Let ${\bf D}$ be a one-dimensional \'{e}tale $(\varphi^r,\Gamma)$-module over
$k_{{\mathcal E}}$. There exists a basis element $g$ for ${\bf D}$, uniquely determined
integers $0\le s({\bf D})\le p-2$ and $1\le n({\bf D})\le p^r-1$ and a
uniquely determined scalar $\xi({\bf D})\in k^{\times}$ such
that \begin{align}\varphi_{{\bf D}}^rg&=\xi({\bf
    D})t^{n({\bf D})+1-p^{r}}g\notag\\\gamma(x)g&-x^{s({\bf
      D})}g\quad\quad\quad\quad\quad\in t\cdot k^+_{{\mathcal E}}\cdot g\notag\end{align}for all
$x\in{\mathbb Z}_p^{\times}$. Thus, one may define $0\le k_{i}({\bf D})\le p-1$ by $n({\bf
  D})=\sum_{i=0}^{r-1}k_{i}({\bf D})p^i$. One has $n({\bf D})\equiv 0$ modulo $(p-1)$.

(b) For any given integers $0\le s\le p-2$ and $1\le n\le p^r-1$ with $n\equiv 0$ modulo $(p-1)$
and any scalar $\xi\in k^{\times}$ there is a uniquely determined
(up to isomorphism) one-dimensional \'{e}tale $(\varphi^r,\Gamma)$-module ${\bf
  D}$ over
$k_{{\mathcal E}}$ with $s=s({\bf D})$ and $n=n({\bf D})$ and $\xi=\xi({\bf D})$.  
\end{lem}

{\sc Proof:} (a) Begin with an arbitrary basis element $g_0$ for ${\bf D}$; then
$\varphi_{{\bf D}}^rg_0=F g_0$ for some unit $F\in k((t))=k_{{\mathcal E}}$. After multiplying $g_0$ with a
suitable power of $t$ we may assume $F=\xi t^{m}(1+t^{n_0}F_0)$ for some
$0\ge m\ge 2-p^r$, some $\xi\in k^{\times}$, some $n_0>0$, and some $F_0\in k[[t]]$ (use
$t^{p^r}\varphi^r=\varphi^rt$). For $g_1=(1+t^{n_0}F_0)g_0$ we then get
$\varphi_{{\bf D}}^rg_1=\xi t^{m}(1+t^{n_1}F_1)g_1$ for some $n_1>n_0>0$ and some $F_1\in
k[[t]]$. We put $g_2=(1+t^{n_1}F_1)g_1$. We may continue in this way. In the
limit we get, by completeness, $g=g_{\infty}\in{\bf D}$ such
that $\varphi_{{\bf D}}^rg=\xi t^{m}g$. It is clear that $\xi({\bf D})=\xi$ and
$n({\bf D})=p^{r}-1+m$ are well defined.

Next, as $\gamma(x)^{p-1}$ (for $x\in{\mathbb Z}_p^{\times}$) is
topologically nilpotent (in the topological group $\Gamma$) and acts by an automorphism on ${\bf D}$, there is some unit $F_x\in
k[[t]]$ with $\gamma(x)g=F_xg$. Thus, the action of $\Gamma$ on ${\bf D}$ induces an action of $\Gamma$ on $k[[t]]/tk[[t]]\cong k$, and this action is given by a homomorphism ${\mathbb Z}_p^{\times}\cong\Gamma\to k^{\times}$, $x\mapsto x^{s({\bf D})}$ for some $0\le s({\bf D})\le p-2$. Computing modulo $t\cdot k^+_{{\mathcal E}}\cdot g$ we have (with $n=n({\bf D})$, $s=s({\bf D})$, $\xi=\xi({\bf D})$)$$x^sg\equiv\gamma(x)g=\xi^{-1}\gamma(x)t^{p^r-n-1}\varphi_{{\bf D}}^rg=\xi^{-1}(\gamma(x)t^{p^r-n-1}\gamma(x^{-1}))\gamma(x)\varphi_{{\bf D}}^rg$$$$\stackrel{(i)}{\equiv}\xi^{-1}(xt)^{p^r-n-1}\gamma(x)\varphi_{{\bf D}}^rg\equiv x^s\xi^{-1}(xt)^{p^r-n-1}\varphi_{{\bf D}}^rg=x^{s+p^r-n-1}g.$$Here $(i)$ follows from $\gamma(x)t^{p^r-n-1}\gamma(x^{-1})-(xt)^{p^r-n-1}\in t^{p^r-n}k[[t]]$, formula (\ref{ennoabschied}). In short, we have $g\equiv x^{p^r-n-1}g$ for all $x\in{\mathbb F}_p^{\times}$. We deduce $n\equiv0$ modulo $p-1$.

(b) Put $D=k[[t]]=k^+_{\mathcal E}$ and $D^*={\rm Hom}_k^{\rm
  ct}(k^+_{\mathcal E},k)$. Endow $D^*$ with an action of $k^+_{\mathcal
  E}$ by putting $(\alpha\cdot\ell)(x)=\ell(\alpha\cdot x)$ for $\alpha\in k^+_{\mathcal
  E}$, $\ell\in D^*$ and $x\in D$. For $j\ge0$ define $\ell_j\in D^*$ by
$$\ell_j(\sum_{i\ge0}a_it^i)=a_j.$$Then $\{\ell_j\}_{j\ge0}$ is a $k$-basis of $D^*$, and we have $t\ell_{j+1}=\ell_j$ and $t\ell_0=0$. Define $\varphi^r:D^*\to D^*$ to be the $k$-linear map with $\varphi^r(\ell_j)=\xi^{-1}\ell_{p^rj+n}$. We then find $$(t^{p^r}(\varphi^r(\ell_j)))(\sum_{i\ge0}a_it^i)=(\xi^{-1}\ell_{p^rj+n})(\sum_{i\ge0}a_it^{i+p^r})=\xi^{-1}a_{p^r(j-1)+n},$$$$(\varphi^r(t\ell_j))(\sum_{i\ge0}a_it^i)=\xi^{-1}\ell_{p^r(j-1)+n}(\sum_{i\ge0}a_it^i)=\xi^{-1}a_{p^r(j-1)+n}$$(with $a_i=0$ resp. $\ell_i=0$ for $i<0$). As $t^{p^r}\varphi^r=\varphi^rt$ in $k^+_{\mathcal
  E}[\varphi^r]$ this means that we have defined an action of $k^+_{\mathcal
  E}[\varphi^r]$ on $D^*$. Next, for $x\in {\mathbb F}_p^{\times}$ and $b,m\in{\mathbb Z}_{\ge0}$ we put $$\gamma(x)\cdot\ell_{-b+n\sum_{i=0}^{m-1}p^{ri}}=x^{-s}(\gamma(x)t^b\gamma(x^{-1}))\ell_{-b+n\sum_{i=0}^{m-1}p^{ri}},$$or equivalently,\begin{gather}\gamma(x)\cdot t^b(\varphi^r)^m\ell_0=x^{-s}(\gamma(x)t^b\gamma(x^{-1}))(\varphi^r)^m\ell_0.\label{enn}\end{gather}We claim that this defines an action of $\gamma$ on $D^*$. In order to check well definedness, i.e. independence on the choice of $b$ and $m$ when only $-b+n\sum_{i=0}^{m-1}p^{ri}$ is given, it is enough to consider the equations $t^b(\varphi^r)^m\ell_0=\xi t^{b+np^{rm}}(\varphi^r)^{m+1}\ell_0$ and to compare the respective effect of $\gamma(x)$ (through formula (\ref{enn})) on either side. Namely, we need to see$$\xi^{-1}x^{-s}(\gamma(x)t^b\gamma(x^{-1}))(\varphi^r)^m\ell_0=x^{-s}(\gamma(x)t^{b+np^{rm}}\gamma(x^{-1}))(\varphi^r)^{m+1}\ell_0,$$or equivalently,$$(\gamma(x)t^b\gamma(x^{-1}))(\xi^{-1}-(\gamma(x)t^{np^{rm}}\gamma(x^{-1}))\varphi^r)(\varphi^r)^m\ell_0=0.$$For this it is enough to see that $\xi^{-1}-(\gamma(x)t^{np^{rm}}\gamma(x^{-1})\varphi^r)$ annihilates $(\varphi^r)^m\ell_0$. Now $n<p^r$ implies $(n+1)p^{rm}> n\sum_{i=0}^mp^{ri}$ and hence that $t^{(n+1)p^{rm}}$ annihilates $(\varphi^r)^{m+1}\ell_0=\xi^{-m-1}\ell_{n\sum_{i=0}^mp^{ri}}$. Therefore it is enough, by formula (\ref{ennoabschied}), to show that $\xi^{-1}-(xt)^{np^{rm}}\varphi^r$ annihilates $(\varphi^r)^m\ell_0$. But $n\equiv 0$ modulo $(p-1)$ means $x^{np^{rm}}=1$. Or claim follows since $\xi^{-1}(\varphi^r)^m\ell_0=t^{np^{rm}}(\varphi^r)^{m+1}\ell_0$. 

In fact, we have defined an action of $k^+_{\mathcal
  E}[\varphi^r,\Gamma]$ on $D^*$. Namely, it is clear from the definitions that the relations $\gamma(x)\cdot t= (\gamma(x)t\gamma(x^{-1}))\cdot \gamma(x)$ in $k^+_{\mathcal
  E}[\varphi^r,\Gamma]$ are respected. That the actions of $\gamma(x)$ and $\varphi^r$ on $D^*$ commute follows easily from the relations $\gamma(x)\varphi^r=\varphi^r\gamma(x)$ and $t^{p^r}\varphi^r=\varphi^rt$ in $k^+_{\mathcal
  E}[\varphi^r,\Gamma]$. Notice that $$\gamma(x)\ell_0=x^{-s}\ell_0\quad\mbox{
    for all }x\in{\mathbb Z}_p^{\times}.$$Now passing to the dual $D\cong (D^*)^*$ of
$D^*$ yields a non degenerate $(\psi^r,\Gamma)$-module over $k^+_{\mathcal E}$ with an
associated \'{e}tale $(\varphi^r,\Gamma)$-module ${\bf
  D}$ over $k_{{\mathcal E}}$ with $s=s({\bf D})$ and $n=n({\bf D})$ and
$\xi=\xi({\bf D})$; this is explained in \cite{dfun} Lemma 6.4. This
dualization argument also proves the uniqueness of ${\bf D}$.\hfill$\Box$\\

{\bf Definition:} We say that an \'{e}tale
$(\varphi^r,\Gamma)$-module ${\bf D}$ over
$k_{{\mathcal E}}$ is $C$-symmetric if it admits a
direct sum decomposition ${\bf D}={\bf D}_1\oplus{\bf D}_2$ with one-dimensional \'{e}tale
$(\varphi^r,\Gamma)$-modules ${\bf D}_1$, ${\bf D}_2$ satisfying the following conditions (1), (2C) and (3C): 

(1) $k_{i}({\bf D}_1)=k_{r-1-i}({\bf D}_2)$ for all $0\le i\le r-1$

(2C) $\xi({\bf D}_1)=\xi({\bf D}_2)$

(3C) $s({\bf D}_2)-s({\bf D}_1)\equiv\sum_{i=0}^{r-1}ik_{i}({\bf D}_1)$ modulo $(p-1)$

(4) $k_{\bullet}({\bf D}_1)\notin\{(0,\ldots,0), (p-1,\ldots,p-1)\}$\\

{\bf Definition:} We say that an \'{e}tale
$(\varphi^r,\Gamma)$-module ${\bf D}$ over
$k_{{\mathcal E}}$ is $B$-symmetric if $r$ is odd and if ${\bf D}$ admits a
direct sum decomposition ${\bf D}={\bf D}_1\oplus{\bf D}_2$ with one-dimensional \'{e}tale
$(\varphi^r,\Gamma)$-modules ${\bf D}_1$, ${\bf D}_2$ satisfying the following conditions (1), (2B) and (3B): 

(1) $k_{i}({\bf D}_1)=k_{r-1-i}({\bf D}_2)$ for all $0\le i\le r-1$, and $k_{\frac{r-1}{2}}({\bf D}_1)=k_{\frac{r-1}{2}}({\bf D}_2)$ is even

(2B) For both ${\bf D}={\bf D}_1$ and ${\bf D}={\bf D}_2$ we
have $\xi({\bf D})=\prod_{i=0}^{r-1}(k_i({\bf D})!)^{-1}$ and $k_i({\bf
  D})=k_{r-1-i}({\bf D})$ for all $1\le i\le \frac{r-1}{2}$

(3B) $s({\bf D}_2)-s({\bf D}_1)\equiv k_{0}({\bf D}_1)-k_{r-1}({\bf D}_1)$ modulo
$(p-1)$

(4) $k_{\bullet}({\bf D}_1)\notin\{(0,\ldots,0), (p-1,\ldots,p-1)\}$\\

\begin{lem} The conjunction of the conditions (1), (2C) and (3C)
(resp. (1), (2B) and (3B)) is symmetric in ${\bf D}_1$ and ${\bf D}_2$. 
\end{lem}

{\sc Proof:} That each one of the conditions (1), (2C) and (2B) is symmetric even individually is obvious. Now $n({\bf
  D})\equiv0$ modulo $p-1$ implies
$\sum_{i=0}^{r-1}k_{i}({\bf D}_1)\equiv0$ modulo $p-1$. Therefore $s({\bf D}_2)-s({\bf
  D}_1)\equiv\sum_{i=0}^{r-1}ik_{i}({\bf D}_1)$ and $k_{i}({\bf D}_1)=k_{r-1-i}({\bf D}_2)$ for
all $i$ (condition (1)) together imply $s({\bf D}_1)-s({\bf
  D}_2)\equiv\sum_{i=0}^{r-1}ik_{i}({\bf D}_2)$. Thus condition (3C) is
symmetric, assuming condition (1). Similarly, condition (3B) is symmetric,
assuming condition (1). \hfill$\Box$\\

{\bf Definition:} (i) Let $\widetilde{\mathfrak S}_C(r)$ denote the set of triples
$(n,s,\xi)$ with integers $1\le n\le p^r-2$ and $0\le s\le p-2$ and scalars
$\xi\in k^{\times}$ such that $n\equiv 0$ modulo $(p-1)$. Let
${\mathfrak S}_C(r)$ denote the quotient of $\widetilde{\mathfrak S}_C(r)$ by
the involution $$(\sum_{i=0}^{r-1}k_ip^i,s,\xi)\mapsto
(\sum_{i=0}^{r-1}k_{r-i-1}p^i,s+\sum_{i=0}^{r-1}ik_i,\xi).$$(Here
and in the following,
in the second component we mean the representative modulo $p-1$ belonging to
$[0,p-2]$.) 

(ii) Let $r$ be odd and
let $\widetilde{\mathfrak S}_B(r)$ denote the set of pairs
$(n,s)$ with integers $1\le n=\sum_{i=0}^{r-1}k_ip^i\le p^r-2$ and $0\le s\le p-2$ such that $n\equiv 0$ modulo $(p-1)$ and such that $k_i=k_{r-1-i}$ for all $1\le i\le
\frac{r-1}{2}$. Let ${\mathfrak S}_B(r)$ denote the quotient of $\widetilde{\mathfrak S}_B(r)$ by
the involution $$(\sum_{i=0}^{r-1}k_ip^i,s)\mapsto
(\sum_{i=0}^{r-1}k_{r-i-1}p^i,s+k_0-k_{r-1}).$$  

\begin{lem}\label{combiflip} (i) Sending ${\bf D}={\bf D}_1\oplus{\bf D}_2$ to $(n({\bf D}_1),
  s({\bf D}_1),\xi({\bf D}_1))$ induces a bijection between the set of
  isomorphism classes of $C$-symmetric \'{e}tale
$(\varphi^r,\Gamma)$-modules and ${\mathfrak S}_C(r)$.

(ii) Sending ${\bf D}={\bf D}_1\oplus{\bf D}_2$ to $(n({\bf D}_1),
  s({\bf D}_1))$ induces a bijection between the set of
  isomorphism classes of $B$-symmetric \'{e}tale
$(\varphi^r,\Gamma)$-modules and ${\mathfrak S}_B(r)$.
\end{lem}

{\sc Proof:} This follows from Lemma \ref{flippa}. \hfill$\Box$\\

{\bf Definition:} Let $r$ be even. Let $\widetilde{\mathfrak S}_D(r)$ denote the set of triples
$(n,s,\xi)$ with integers $1\le n=\sum_{i=0}^{r-1}k_ip^i\le p^r-2$ and $0\le s\le p-2$ and scalars
$\xi\in k^{\times}$ such that $n\equiv 0$ modulo $(p-1)$ and such that
$k_i=k_{i+\frac{r}{2}}$ for all $1\le i\le \frac{r}{2}-2$. We consider the following permutations $\iota_0$ and
$\iota_1$ of $\widetilde{\mathfrak S}_D(r)$. The value of $\iota_0$ at
$(\sum_{i=0}^{r-1}k_ip^i,s,\xi)$
is$$(k_{\frac{r}{2}}+\sum_{i=1}^{\frac{r}{2}-2}k_ip^i+k_{r-1}p^{\frac{r}{2}-1}+k_0p^{\frac{r}{2}}+\sum_{i=\frac{r}{2}+1}^{r-2}k_ip^i+k_{\frac{r}{2}-1}p^{r-1},s+\sum_{i=0}^{\frac{r}{2}-1}k_i,\xi).$$The value of $\iota_1$ at
$(\sum_{i=0}^{r-1}k_ip^i,s,\xi)$
is$$(\sum_{i=0}^{r-1}k_{r-i-1}p^i,s+\frac{r-2}{4}(k_{\frac{r}{2}}+k_0)+\sum_{i=2}^{\frac{r}{2}-1}(i-1)k_{\frac{r}{2}-i},\xi)$$if
$\frac{r}{2}$ is odd, whereas if $\frac{r}{2}$ is even the value
is $$(k_{\frac{r}{2}-1}+\sum_{i=1}^{\frac{r}{2}-2}k_{r-i-1}p^i+k_{\frac{r}{2}}p^{\frac{r}{2}-1}+k_{r-1}p^{\frac{r}{2}}+\sum_{i=\frac{r}{2}+1}^{r-2}k_{r-i-1}p^i+k_0p^{r-1},s+(\frac{r}{4}-1)k_{\frac{r}{2}}+\frac{r}{4}k_0+\sum_{i=2}^{\frac{r}{2}-1}(i-1)k_{\frac{r}{2}-i}p^i,\xi).$$It
is straightforward to check that $\iota_0^2={\rm id}$ and $\iota_0\iota_1=\iota_1\iota_0$, and moreover that
$\iota_1^2={\rm id}$ if $\frac{r}{2}$ is odd, but
$\iota_1^2=\iota_0$ if $\frac{r}{2}$ is even. In either case, the subgroup $\langle\iota_0,\iota_1\rangle$ of ${\rm
  Aut}({\widetilde{\mathfrak S}_D(r)})$ generated by $\iota_0$ and $\iota_1$ is
commutative and contains $4$ elements. We let ${\mathfrak S}_D(r)$ denote the
quotient of $\widetilde{\mathfrak S}_D(r)$ by the action of $\langle\iota_0,\iota_1\rangle$.\\
 
{\bf Definition:} Let $r$ be even. We say that an \'{e}tale
$(\varphi^r,\Gamma)$-module ${\bf D}$ over
$k_{{\mathcal E}}$ is $D$-symmetric if it admits a
direct sum decomposition ${\bf D}={\bf D}_{11}\oplus{\bf D}_{12}\oplus{\bf D}_{21}\oplus{\bf D}_{22}$ with one-dimensional \'{e}tale
$(\varphi^r,\Gamma)$-modules ${\bf D}_{11}$, ${\bf D}_{12}$, ${\bf
  D}_{21}$, ${\bf D}_{22}$ satisfying the following conditions:

(1) For all $1\le i\le \frac{r}{2}-2$ and all $1\le s, t\le 2$ we have $k_i({\bf D}_{st})=k_{\frac{r}{2}+i}({\bf D}_{st})$

(2) For all $1\le i\le \frac{r}{2}-2$ we have $k_i({\bf D}_{11})=k_i({\bf D}_{12})$ and
$k_i({\bf D}_{21})=k_i({\bf D}_{22})$

(3) For $j=1$ and $j=2$ we have $$k_0({\bf D}_{j1})=k_{\frac{r}{2}}({\bf D}_{j2}),\quad k_{\frac{r}{2}}({\bf D}_{j1})=k_0({\bf D}_{j2}),\quad k_{\frac{r}{2}-1}({\bf D}_{j1})=k_{r-1}({\bf D}_{j2}),\quad k_{r-1}({\bf D}_{j1})=k_{\frac{r}{2}-1}({\bf D}_{j2}).$$ 

(4) $$k_i({\bf D}_{11})=k_{r-i-1}({\bf
  D}_{21})\quad\mbox{ and }\quad k_i({\bf D}_{12})=k_{r-i-1}({\bf
  D}_{22})$$if $i\in[0,r-1]$ and $\frac{r}{2}$ is odd, or if
$i\in[1,\frac{r}{2}-2]\cup [\frac{r}{2}+1,r-2]$ and $\frac{r}{2}$ is
even. Moreover, if $\frac{r}{2}$ is
even then$$k_0({\bf D}_{11})=k_{r-1}({\bf D}_{21}),\quad
k_{\frac{r}{2}-1}({\bf D}_{11})=k_{0}({\bf D}_{21}),\quad k_{\frac{r}{2}}({\bf
  D}_{11})=k_{\frac{r}{2}-1}({\bf D}_{21}),\quad k_{r-1}({\bf
  D}_{11})=k_{\frac{r}{2}}({\bf D}_{21})$$

(5) $\xi({\bf D}_{11})=\xi({\bf D}_{12})=\xi({\bf D}_{21})=\xi({\bf D}_{22})$

(6) Modulo $(p-1)$ we have\begin{gather}s({\bf D}_{j2})-s({\bf D}_{j1})\equiv\sum_{i=0}^{\frac{r}{2}-1}k_i({\bf D}_{j1})\quad\quad\mbox{ for }j=1,2\notag\\s({\bf D}_{21})-s({\bf D}_{11})\equiv\left\{\begin{array}{l@{\quad:\quad}l}\frac{r-2}{4}(k_{\frac{r}{2}}({\bf
  D}_{11})+k_0({\bf
  D}_{11}))+\sum_{i=2}^{\frac{r}{2}-1}(i-1)k_{\frac{r}{2}-i}({\bf
  D}_{11})&\frac{r}{2}\mbox{ is odd}\notag\\(\frac{r}{4}-1)k_{\frac{r}{2}}({\bf
  D}_{11})+\frac{r}{4}k_0({\bf
  D}_{11})+\sum_{i=2}^{\frac{r}{2}-1}(i-1)k_{\frac{r}{2}-i}({\bf
  D}_{11})&\frac{r}{2}\mbox{ is even}\end{array}\right.\notag\end{gather}

(7) $k_{\bullet}({\bf D}_{11})\notin\{(0,\ldots,0), (p-1,\ldots,p-1)\}$\\

\begin{lem}\label{dflip} Sending ${\bf D}={\bf D}_{11}\oplus{\bf D}_{12}\oplus{\bf D}_{21}\oplus{\bf D}_{22}$ to $(n({\bf D}_{11}),
  s({\bf D}_{11}),\xi({\bf D}_{11}))$ induces a bijection between the set of
  isomorphism classes of $D$-symmetric \'{e}tale
$(\varphi^r,\Gamma)$-modules and ${\mathfrak S}_D(r)$.
\end{lem}

{\sc Proof:} Again we use Lemma \ref{flippa}. For a one-dimensional \'{e}tale $(\varphi^r,\Gamma)$-module
${\bf D}$ over $k_{{\mathcal E}}$ put $\alpha({\bf D})=(n({\bf D}),
  s({\bf D}),\xi({\bf D}))$. If ${\bf D}={\bf D}_{11}\oplus{\bf D}_{12}\oplus{\bf
    D}_{21}\oplus{\bf D}_{22}$ is $D$-symmetric as above, then $\alpha({\bf D}_{st})$ is an element of $\widetilde{\mathfrak
    S}_D(r)$, for all $1\le s,t\le2$. Moreover, it is
  straightforward to check $\iota_0(\alpha({\bf D}_{11}))=\alpha({\bf
    D}_{12})$, $\iota_0(\alpha({\bf D}_{21}))=\alpha({\bf
    D}_{22})$, $\iota_1(\alpha({\bf D}_{11}))=\alpha({\bf
    D}_{21})$ and $\iota_1(\alpha({\bf D}_{12}))=\alpha({\bf
    D}_{22})$. It follows that the above map is well defined and bijective.\hfill$\Box$\\

{\bf Definition:} We say that an \'{e}tale
$(\varphi^r,\Gamma)$-module ${\bf D}$ over
$k_{{\mathcal E}}$ is $A$-symmetric if ${\bf D}$ admits a
direct sum decomposition ${\bf D}=\oplus_{i=0}^{r-1}{\bf D}_i$ with
one-dimensional \'{e}tale $(\varphi^r,\Gamma)$-modules ${\bf D}_i$ satisfying$$k_i({\bf D}_j)=k_{i-j}({\bf D}_0),\quad\quad \xi({\bf
  D}_j)=\xi({\bf D}_0),\quad\quad s({\bf D}_0)-s({\bf
  D}_j)\equiv\sum_{i=1}^jk_{-i}({\bf D}_0)\mbox{ modulo }(p-1)$$ for all $i$, $j$ (where we understand the sub index in $k_{?}$ as
the unique representative in $[0,r-1]$ modulo $r$), and moreover $$k_{\bullet}({\bf D}_0)\notin\{(0,\ldots,0), (p-1,\ldots,p-1)\}.$$

\section{Semiinfinite chamber galleries and the functor ${\bf D}$}

\subsection{Power multiplicative elements in the extended affine Weyl group}

\label{powbas}

Let $G$ be the group of ${\mathbb Q}_p$-rational points of a ${\mathbb
  Q}_p$-split connected reductive group over ${\mathbb Q}_p$. Fix a maximal
${\mathbb Q}_p$-split torus $T$ in $G$, let $N(T)$ be its normalizer in $G$. Let $\Phi$ denote the set of roots of $T$. For
$\alpha\in\Phi$ let $N_{\alpha}$ be the corresponding root subgroup in
$G$. Choose a positive system $\Phi^+$ in $\Phi$, let $\Delta\subset\Phi^+$ be
the set of simple roots. Let $N=\prod_{\alpha\in\Phi^+}N_{\alpha}$. 

Let $X$ denote the semi simple Bruhat-Tits building of $G$, let $A$ denote
 its apartment corresponding to $T$. Our notational and terminological convention is that the facets of $A$ or
 $X$ are {\it closed} in $X$ (i.e. {\it contain} all their faces (the lower dimensional facets
 at their boundary)). A chamber is a facet of codimension $0$. For a chamber $D$ in $A$ let $I_D$ be the Iwahori subgroup in $G$ fixing $D$. We notice that $I_D\cap N=\prod_{\alpha\in\Phi^+}I_{D}\cap
N_{\alpha}$.

Fix a special vertex $x_0$ in $A$, let $K$ be the corresponding hyperspecial
maximal compact open subgroup in $G$. Let $T_0=T\cap K$ and $N_0=N\cap K$. We have the isomorphism
$T/T_0\cong X_*(T)$ sending $\xi\in X_*(T)$ to the class of
$\xi(p)\in T$. Let $I\subset K$ be the Iwahori subgroup determined by $\Phi^+$. [If ${\rm
  red}:K\to\overline{K}$ denotes the reduction map onto the reductive (over
${\mathbb F}_p)$ quotient $\overline{K}$ of $K$, then $I={\rm red}^{-1}({\rm
  red}(T_0N_0))$.] Let $C\subset A$ be the chamber fixed by $I$. Thus $I=I_C$ and $I\cap N=N_0=\prod_{\alpha\in\Phi^+}N_0\cap N_{\alpha}$. 

We are interested in semiinfinite chamber galleries\begin{gather}C^{(0)},
  C^{(1)},C^{(2)},C^{(3)},\ldots\label{cgall}\end{gather}in $A$ such that
$C=C^{(0)}$ (and thus $I=I_{C^{(0)}}$) and such that,
setting $$N_0^{(i)}=I_{C^{(i)}}\cap N,$$we have $N_0=N_0^{(0)}$ and \begin{gather}N_0^{(0)}\supset N_0^{(1)}\supset
  N_0^{(2)}\supset N_0^{(3)}\supset\ldots\quad\mbox{ with
  }[N_0^{(i)}:N_0^{(i+1)}]=p\mbox{ for all }i\ge0.\label{n0ineu}\end{gather}

For any two chambers $D_1\ne D_2$ in $A$ sharing a common facet of codimension $1$ there is a uniquely determined $\delta\in\Phi$ with $I_{D_1}\cap
N_{\delta'}=I_{D_2}\cap
N_{\delta'}$ for $\delta'\notin\{\delta,-\delta\}$, with $I_{D_1}\cap
N_{-\delta}\subset I_{D_2}\cap
N_{-\delta}$ and with $I_{D_2}\cap
N_{\delta}\subset I_{D_1}\cap
N_{\delta}$, and both these inclusions are of index $p$. Moreover, the one-codimensional facet shared by $D_1$ and $D_2$ is contained in a translate of the hyperplane corresponding to $\delta$ (or equivalently, corresponding to $-\delta$).\footnote{Pick $g\in N(T)$ acting as the reflection in the affine hyperplane which contains the one-codimensional facet shared by $D_1$ and $D_2$. Then $g D_1=D_2$ and hence $I_{D_2}=g I_{D_1}g^{-1}$, so we can read off all these claims.}

Applying this remark to $D_1=C^{(i)}$, $D_2=C^{(i+1)}$ and putting $\alpha^{(i)}=\delta$ we see that condition (\ref{n0ineu}) precisely means that $\alpha^{(i)}\in\Phi^+$ for all $i\ge0$. Conversely, if condition (\ref{n0ineu}) holds true then the sequence
$\alpha^{(0)},\alpha^{(1)},\alpha^{(2)},\alpha^{(3)},\ldots$ in $\Phi^+$ can be characterized as follows: Setting$$e[i,{\alpha}]=|\{0\le j\le i-1\,|\,\alpha=\alpha^{(j)}\}|$$for
$i\ge0$ and $\alpha\in\Phi^+$, we
have\begin{gather}N_0^{(i)}=\prod_{\alpha\in\Phi^+}(N_0\cap
  N_{\alpha})^{p^{e[i,\alpha]}}.\notag\end{gather}

Suppose that the center $Z$ of $G$ is connected. Then $G/Z$ is a semisimple
group of adjoint type with maximal torus $\check{T}=T/Z$. Let
$\check{T}_0=T_0/(T_0\cap Z)\subset \check{T}$. The extended affine Weyl group $\widehat{W}=N(\check{T})/\check{T}_0$ can be identified with the
semidirect product between the finite Weyl group
$W=N(\check{T})/\check{T}=N(T)/T$ and $X_*(\check{T})$. We identify $A=
X_*(\check{T})\otimes{\mathbb R}$ such that $x_0\in A$ corresponds to the
origin in the ${\mathbb R}$-vector space $X_*(\check{T})\otimes{\mathbb
  R}$. We then regard $\widehat{W}$ as acting on $A$ through affine
transformations. We regard $\Delta\subset X^*({T})$ as a subset of
$X^*(\check{T})$. We usually enumerate the elements of $\Delta$ as $\alpha_1,\ldots,\alpha_d$, and we
enumerate the corresponding simple reflection $s_{\alpha}\in W$ for
$\alpha\in\Delta$ as $s_1,\ldots,s_d$ with $s_i=s_{\alpha_i}$. We assume that the root system $\Phi$ is
irreducible. Let ${\alpha}_0\in\Phi$ be the {\it negative} of the highest
root. Let $s_{\alpha_0}$ be the corresponding reflection in the finite Weyl
group $W$; define the affine reflection $s_0=t_{\alpha_0^{\vee}}\circ s_{\alpha_0}\in
\widehat{W}$, where $t_{\alpha_0^{\vee}}$ denotes the translation by the
coroot $\alpha_0^{\vee}\in A$ of $\alpha_0$. The affine Weyl group
$W_{{\rm aff}}$ is the subgroup of $\widehat{W}$ generated by
$s_0,s_1,\ldots,s_d$; in fact it is a Coxeter group with these Coxeter generators. The corresponding length function $\ell$ on $W_{{\rm
    aff}}$ extends to $\widehat{W}$. 

The above discussion shows that if the gallery (\ref{cgall}) satisfies condition (\ref{n0ineu}) and if $w\in \widehat{W}$ satisfies $C^{(i)}=wC$ for some $i\ge0$, then $$p^{\ell(w)}=p^i=p^{[N_0:I_{C^{(i)}}\cap N_0]}.$$

 Let $X_*(\check{T})_+$ denote the set of dominant cocharacters.  [Let $T_+=\{t\in T\,|\,t N_0
t^{-1}\subset N_0\}$, then $X_*(\check{T})_+$ is the image of $T_+$ under the map
$T_+\subset T\to T/T_0\cong X_*(T)\to X_*(\check{T})$.] The
monoid $X_*(\check{T})_+$ is free and has a distinguished basis $\nabla$, the set of
fundamental coweights. The cone ${\mathcal C}$
(vector chamber) in $A$ with origin in $x_0$ which is spanned by all the
$-\xi$ for $\xi\in \nabla$ contains $C$, and $C$ is precisely the 'top'
chamber of this cone. The reflections $s_0,s_1,\ldots,s_d$ are precisely the
reflections in the affine hyperplanes (walls) of $A$ which contain a
codimension-$1$-face of $C$.\\

%
%
%

{\bf Definition:} We say that $w\in \widehat{W}$ is power multiplicative if
$\ell(w^m)=m\cdot\ell(w)$ for all $m\ge0$.\footnote{After finishing this article we found out that power multiplicative elements of $\widehat{W}$ have been studied e.g. in \cite{hn} under the name of {\it straight} elements.}\\

Of course, any element in the image
of $T\to N(\check{T})\to \widehat{W}=N(\check{T})/\check{T}_0$ is
power multiplicative because it acts on $A$ by translation (and $\ell(w)$ is the gallery distance between $C$ and $wC$). 

Suppose we are given a fundamental coweight $\tau\in \nabla$ and some non trivial element $\phi\in
\widehat{W}$ satisfying the following conditions:

(a) $\phi$ is power multiplicative,

(b) $\tau$ is minuscule, i.e. we have $\langle\alpha,\tau\rangle\in\{0,1\}$
for all $\alpha\in\Phi^+$,

(c) viewing $\tau$ via the embedding $X_*(\check{T})\subset\widehat{W}$ as an element
of $\widehat{W}$, we
have\begin{gather}\phi^{\mathbb
    N}\cap \tau^{\mathbb N}\ne \emptyset.\label{potenzidend}\end{gather}

\begin{lem}\label{concept} Let $\phi$ and $\tau$ be as above. Write
  $\phi=\phi'v$ with $\phi'\in W_{{\rm aff}}$ and $v\in \widehat{W}$ with
  $vC=C$. Choose a reduced
  expresssion$$\phi'={s}_{\beta(1)}\cdots{s}_{\beta(r)}$$of $\phi'$ with some function
  $\beta:\{1,\ldots,r\}\to\{0,\ldots,d\}$ (with $r=\ell(\phi)=\ell(\phi')$) and
  put\begin{gather}C^{(ar+b)}=\phi^as_{\beta(1)}\cdots s_{\beta(b)}C\label{neuihg}\end{gather}for $a,b\in{\mathbb
    Z}_{\ge0}$ with $0\le b<r$. Lift $\tau\in \nabla\subset X_*(\check{T})$
  to some element of $X_*({T})$ and denote again by $\tau$ the
  corresponding homomorphism ${\mathbb Z}_p^{\times}\to T_0$. Then we have:

(i) The sequence$$C=C^{(0)}, C^{(1)},
C^{(2)},\ldots$$satisfies hypothesis (\ref{n0ineu}). In particular we may define $\alpha^{(j)}\in\Phi^+$
for all $j\ge0$. 

(ii) For any $j\ge0$ we have $\alpha^{(j)}\circ\tau={\rm id}_{{\mathbb
  Z}_p^{\times}}$. 

(iii) For any lifting $\phi\in N(T)$ of $\phi\in
\widehat{W}$ we have $\tau(a)\phi=\phi\tau(a)$ in $N(T)$, for all
$a\in{\mathbb Z}_p^{\times}$.
\end{lem}

{\sc Proof:} Given a chamber gallery (\ref{cgall}) in $A$ with $C=C^{(0)}$ and $N_0=N_0^{(0)}$, condition (\ref{n0ineu}) is equivalent with saying that (\ref{cgall}) is a minimal chamber gallery in the cone ${\mathcal C}$, i.e. all $C^{(i)}$ are contained in ${\mathcal C}$, and the gallery distance from $C=C^{(0)}$ to $C^{(i)}$ is $i$. Conversely, notice that if $C^{(i)}$ belongs to ${\mathcal C}$ and the gallery distance from $C=C^{(0)}$ to $C^{(i)}$ is $i$, then also all $C^{(j)}$ for $0\le j\le i$ belong to ${\mathcal C}$. (Otherwise the gallery would have to cross one of the bounding hyperplanes of ${\mathcal C}$ and then cross it again later in the reverse direction; but a minimal gallery never a crosses a wall (hyperplane translate) back and forth.) Recall that, assuming condition (\ref{n0ineu}), the set $\{\alpha^{(j)}\,|\,j\ge0\}$ is precisely the set of all $\alpha\in\Phi^+$ such that the gallery (\ref{cgall}) crosses a translate of the hyperplane corresponding to $\alpha$. If there is some $\epsilon\in X_*(T)$ and some $s\in{\mathbb N}$ such that $C^{(i+s)}=\epsilon C^{(i)}$ for all $i$, then, since the action of $\epsilon$ is by translation on $A$, the latter set is the same as the set of all $\alpha\in\Phi^+$ such that the gallery $C=C^{(0)},C^{(1)},\ldots, C^{(s)}$ crosses a translate of the hyperplane corresponding to $\alpha$. 

Now let (\ref{cgall}) be the chamber gallery defined by formula (\ref{neuihg}).  
(i) If $n\in{\mathbb N}$ is such that $\phi^n\in X_*(T)$, then the gallery distance between $C$ and $\phi^nC$ is $\ell(\phi^n)$. By construction, the gallery $C=C^{(0)},C^{(1)},\ldots, C^{(n\ell(\phi))}=\phi^nC$ has length $n\ell(\phi)$, i.e. length $\ell(\phi^n)$ as $\phi$ is power multiplicative. Hence it must be a minimal chamber gallery. As $\tau$ is minuscule, it is in particular a dominant
coweight. Therefore it follows from hypothesis (\ref{potenzidend}) that also
some power $\phi^n$ of $\phi$ is a dominant
coweight. We thus obtain that $C^{(n'\ell(\phi))}=\phi^{n'}C$ lies in ${\mathcal C}$ for each $n'\in{\mathbb N}$ divisible by $n$, but then, as the gallery is minimal and ${\mathcal C}$ is a cone, $C^{(i)}$ lies in ${\mathcal C}$ for each $i\ge0$. We get statement (i) in view of the preceeding discussion.

(ii) The preceeding discussion shows $\{\alpha^{(j)}\,|\,j\ge0\}=\{\alpha\in\Phi^+\,|\,\langle\alpha,\phi^m\rangle\ne0\}$ for any $m\in{\mathbb N}$ such that $\phi^m$ belongs to $X_*(T)$. As $\phi$ is power multiplicative, hypothesis (\ref{potenzidend}) furthermore gives$$\{\alpha^{(j)}\,|\,j\ge0\}=\{\alpha\in\Phi^+\,|\,\langle\alpha,\phi^m\rangle\ne0\}=\{\alpha\in\Phi^+\,|\,\langle\alpha,\tau\rangle\ne0\}$$and
as $\tau$ is minuscule this is the set $$\{\alpha\in\Phi^+\,|\,\langle\alpha,\tau\rangle=1\}=\{\alpha\in\Phi^+\,|\,\alpha\circ\tau={\rm id}_{{\mathbb Z}_p^{\times}}\}.$$

(iii) By hypothesis (\ref{potenzidend}) we have $\tau^m=\phi^n$ for some
$m,n\in{\mathbb N}$. We deduce
$\tau^m=\phi\tau^m\phi^{-1}=(\phi\tau\phi^{-1})^{m}$ and hence also $\tau=\phi\tau\phi^{-1}$ as $\tau$ and
$\phi\tau\phi^{-1}$ belong to the free abelian group $X_*(T)$. Thus
$\tau\phi=\phi\tau$ in $\widehat{W}$ which implies claim (iii).\hfill$\Box$\\

{\bf Remark:} For a given minuscule fundamental coweight $\tau\in \nabla$ any positive power $\tau^m$ of $\tau$ belongs to
$X_*(\check{T})$, and so $\phi=\tau^m$ satisfies the assumptions of Lemma
\ref{concept}. However, for our purposes it is of interest to find $\phi$ (as
in Lemma \ref{concept}, possibly also required to project to $W_{\rm aff}$) of small length; the minimal positive power of
$\tau$ belonging to $X_*(\check{T})$ is usually not optimal in this sense.

\subsection{The functor ${\bf D}$}

\label{erklaer}

By $I_0$ we denote the pro-$p$-Iwahori subgroup contained in
$I$. We often read $\overline{T}=T_0/T_0\cap I_0$ as a subgroup of $T_0$ by means of the
Teichm\"uller character. Conversely, we read characters of $\overline{T}$ also
as characters of $T_0$ (and do not introduce another name
for these inflations). 

Let ${\rm ind}_{I_0}^{G}{\bf 1}_{{\mathfrak o}}$ denote the ${\mathfrak o}$-module of
${\mathfrak o}$-valued compactly supported functions $f$ on $G$ such that $f(ig)=f(g)$ for all
$g\in G$, all $i\in I_0$. It is a $G$-representation by means of
$(g'f)(g)=f(gg')$ for $g,g'\in G$. Let $${\mathcal
  H}(G,I_0)={\rm End}_{{\mathfrak
    o}[G]}({\rm ind}_{I_0}^{G}{\bf 1}_{{\mathfrak o}})^{\rm op}$$denote the corresponding
pro-$p$-Iwahori Hecke algebra with coefficients in ${\mathfrak o}$. For a subset $H$ of $G$ let $\chi_H$ denote the
characteristic function of $H$. For $g\in G$ let $T_g\in {\mathcal
  H}(G,I_0)$ denote the Hecke operator corresponding to the double coset
$I_0gI_0$. It sends $f:G\to{\mathfrak o}$
to $$T_g(f):G\longrightarrow{\mathfrak o},\quad\quad h\mapsto\sum_{x\in
  I_0\backslash G}\chi_{I_0gI_0}(hx^{-1})f(x).$$Let ${\rm Mod}^{\rm fin}({\mathcal
  H}(G,I_0))$ denote the category of ${\mathcal
  H}(G,I_0)$-modules which as ${\mathfrak o}$-modules are of finite length. We write ${\mathcal H}(G,I_0)_k={\mathcal H}(G,I_0)\otimes_{\mathfrak
  o}k$. Given liftings $\dot{s}\in N(T)$ of all $s\in S$ we let ${\mathcal H}(G,I_0)_{{\rm aff},k}$ denote the $k$-subalgebra of
${\mathcal H}(G,I_0)_k$ generated by the $T_{\dot{s}}$ for all $s\in{S}$ and the $T_t$ for $t\in \overline{T}$.

Suppose we are given a reduced expression\begin{gather}\phi=\epsilon\dot{s}_{\beta(1)}\cdots\dot{s}_{\beta(r)}\label{einlfixred}\end{gather}(some function
  $\beta:\{1,\ldots,r=\ell(\phi)\}\to\{0,\ldots,d\}$, some $\epsilon\in Z$) of a power multiplicative
  element $\phi\in N(T)$, some power of which maps to a
dominant coweight in $N(T)/ZT_0$. Put$$C^{(ar+b)}=\phi^as_{\beta(1)}\cdots s_{\beta(b)}C$$for
$a,b\in{\mathbb Z}_{\ge0}$ with $0\le b<r$. Then, by power multiplicativity of
$\phi$, the sequence (\ref{cgall}) thus defined satisfies property
(\ref{n0ineu}), cf. Lemma \ref{concept}. Therefore we may use it to place ourselves into the setting (and
notations) of \cite{dfun}, as follows. 

We define the half tree $Y$ whose edges are the elements in the
$N_0$-orbits of the $C^{(i)}\cap C^{(i+1)}$ and whose vertices are the elements in the
$N_0$-orbits of the $C^{(i)}$. We choose an isomorphism $\Theta:Y\cong{\mathfrak
  X}_+$ with the $\lfloor{\mathfrak
  N}_0,\varphi,\Gamma\rfloor$-equivariant half sub tree ${\mathfrak X}_+$ of
the Bruhat Tits tree of ${\rm GL}_2({\mathbb Q}_p)$,
satisfying the requirements of Theorem 3.2 of loc.cit.. It sends the edge
$C^{(i)}\cap C^{(i+1)}$ (resp. the vertex $C^{(i)}$) of $Y$ to the edge
${\mathfrak e}_{i+1}$ (resp. the vertex ${\mathfrak v}_i$) of ${\mathfrak
  X}_+$. The half tree $\overline{\mathfrak X}_+$ is obtained from ${\mathfrak X}_+$ by removing the 'loose' edge ${\mathfrak e}_{0}$.

To an ${\mathcal
  H}(G,I_0)$-module $M$ we associate the $G$-equivariant (partial)
coefficient system ${\mathcal V}_M^X$ on $X$. Briefly, its value at the
chamber $C$ is ${\mathcal
  V}_M^X(C)=M$. The transition maps ${\mathcal
  V}_M^X(D)\to{\mathcal
  V}_M^X(F)$ for chambers (codimension-$0$-facets) $D$ and codimension-$1$-facets $F$ with $F\subset D$
are injective, and ${\mathcal
  V}_M^X(F)$ for any such $F$ is the sum of the images of the ${\mathcal
  V}_M^X(D)\to{\mathcal
  V}_M^X(F)$ for all $D$ with $F\subset D$. 

The pushforward $\Theta_*{\mathcal V}_M$ of the restriction of ${\mathcal
  V}_M^X$ to $Y$ is a $\lfloor{\mathfrak
  N}_0,\varphi^r,\Gamma_0\rfloor$-equivariant coefficient system on
${\mathfrak X}_+$ (and by further restriction on $\overline{\mathfrak
  X}_+$). This leads to the exact functor\begin{gather}M\mapsto {\bf
  D}(\Theta_*{\mathcal V}_M)\label{dfual}\\{\bf
  D}(\Theta_*{\mathcal V}_M)=H_0(\overline{\mathfrak
  X}_+,\Theta_*{\mathcal V}_M)\otimes_{{\mathcal O}_{\mathcal E}^+}{\mathcal O}_{\mathcal E}=H_0({\mathfrak
  X}_+,\Theta_*{\mathcal V}_M)\otimes_{{\mathcal O}_{\mathcal E}^+}{\mathcal O}_{\mathcal E}\notag\end{gather}from ${\rm Mod}^{\rm fin}({\mathcal
  H}(G,I_0))$ to the
category of $(\varphi^{r},\Gamma_0)$-modules over ${\mathcal O}_{\mathcal E}$, where $r=\ell(\phi)$. If in addition we are given a homomorphism $\tau:{\mathbb Z}_p^{\times}\to
T_0$ satisfying the conclusions of Lemma \ref{concept} (with respect to
$\phi$), then this functor in fact takes values in the category of
$(\varphi^{r},\Gamma)$-modules over ${\mathcal O}_{\mathcal E}$. 

For $0\le i\le r-1$ we put $$y_i=\dot{s}_{\beta(1)}\cdots\dot{s}_{\beta(i+1)}\dot{s}^{-1}_{\beta(i)}\cdots\dot{s}^{-1}_{\beta(1)}.$$

\begin{lem}\label{prepa} (a) $y_i=y_{i-1}\cdots y_0\dot{s}_{\beta(i+1)}y_0^{-1}\cdots
  y_{i-1}^{-1}$ and $\phi=\epsilon y_{r-1}\cdots y_0$.

(b) $y_i$ is the affine reflection in the wall passing through
$C^{(i)}\cap C^{(i+1)}$. 

\end{lem} 

{\sc Proof:} To see (b) observe that $y_i$ indeed {\it is} a reflection, and
that it sends $C^{(i)}$ to $C^{(i+1)}$. \hfill$\Box$\\

\subsection{Supersingular modules}

\label{basic}

For $\alpha\in\Phi$ we denote by $\alpha^{\vee}$ the associated coroot. For
any $\alpha\in\Phi$ there is a corresponding homomorphism of algebraic groups
$\iota_{\alpha}:{\rm SL}_2({\mathbb Q}_p)\to G$ as described in \cite{jan}, Ch.II,
section 1.3. The element $\iota_{\alpha}(\nu)$ belongs to $I\cap N_{\alpha}$ and
generates it as a topological group. For $x\in{\mathbb F}_p^{\times}\subset
{\mathbb Z}_p^{\times}$ (via the Teichm\"uller character) we
have $\alpha^{\vee}(x)=\iota_{\alpha}(h(x))\in T$.

For $0\le i\le d$ let $\overline{T}_{\alpha^{\vee}_i}$ denote the subgroup of $\overline{T}$ given by the ${\mathbb F}_p$-valued points of the schematic closure of $\alpha^{\vee}_i(\overline{\mathbb F}_p^{\times})$ (in the obvious ${\mathbb F}_p$-group scheme (split torus) underlying $\overline{T}$). For a character $\lambda:\overline{T}\to k^{\times}$ let $S_{\lambda}$ be the subset of $S=\{s_0,\ldots,s_d\}$ consisting of all $s_i$ such that
$\lambda|_{\overline{T}_{\alpha^{\vee}_i}}$ is trivial. Given
$\lambda$ and a
subset ${\mathcal J}$ of $S_{\lambda}$ there
is a uniquely determined character $$\chi_{\lambda,{\mathcal
  J}}:{\mathcal H}(G,I_0)_{{\rm
    aff},k}\longrightarrow k$$which sends $T_t$ to $\lambda(t^{-1})$ for $t\in
\overline{T}$, which sends $T_{\dot{s}}$ to $0$ for $s\in {S}-{\mathcal
  J}$ and which
sends $T_{\dot{s}}$ to $-1$ for $s\in {\mathcal
  J}$ (see \cite{vigneras} Proposition 2). For $0\le i\le d$ we define a number $0\le k_i=k_i(\lambda,{\mathcal
   J})\le p-1$ such that\begin{gather}\lambda(\alpha_i^{\vee}(x))=x^{k_i}\quad\quad\mbox{ for all
}x\in{\mathbb F}_p^{\times},\label{kljfmlc}\end{gather}as follows. If $s_i\in{\mathcal
  J}$ put $k_i=p-1$. Otherwise let $k_i$ be the unique integer in $[0,p-2]$ satisfying formula (\ref{kljfmlc}).\\

Vign\'{e}ras defined the notion of a {\it supersingular} ${\mathcal
  H}(G,I_0)_k$-module, see \cite{vigneras}, \cite{ollcomp}. Let ${\mathcal H}(G,I_0)'_{{\rm
    aff},k}$ denote the $k$-subalgebra of ${\mathcal
  H}(G,I_0)_k$ generated by ${\mathcal H}(G,I_0)_{{\rm
    aff},k}$ together with all $T_z$ for $z\in Z$.

\begin{satz}\label{classsusi} (Vign\'{e}ras, Ollivier)

(a) Let $\chi:{\mathcal H}(G,I_0)'_{{\rm
    aff},k}\to k$ be a character extending $\chi_{\lambda,{\mathcal
  J}}:{\mathcal H}(G,I_0)_{{\rm
    aff},k}\longrightarrow k$ for some $\lambda$, ${\mathcal
  J}$. If $k_{\bullet}\notin\{(0,\ldots,0),(p-1,\ldots,p-1)\}$ then the ${\mathcal
  H}(G,I_0)_k$-module ${\mathcal
  H}(G,I_0)_k\otimes_{{\mathcal H}(G,I_0)'_{{\rm
    aff},k}}\chi$ is supersingular.

(b) Any supersingular ${\mathcal
  H}(G,I_0)_k$-module of finite length contains a character $\chi_{\lambda,{\mathcal
  J}}:{\mathcal H}(G,I_0)_{{\rm
    aff},k}\longrightarrow k$ for some $\lambda$, ${\mathcal
  J}$ with $k_{\bullet}\notin\{(0,\ldots,0),(p-1,\ldots,p-1)\}$.
\end{satz}

{\sc Proof:} See \cite{ollcomp} Theorem 5.14, Corollary 5.16, at least for the case where the root system underlying $G$ is irreducible. See \cite{visusi} for more general statements. \hfill$\Box$\\

It follows that the simple supersingular ${\mathcal
  H}(G,I_0)_k$-modules are precisely the simple quotients of the ${\mathcal
  H}(G,I_0)_k$-modules ${\mathcal
  H}(G,I_0)_k\otimes_{{\mathcal H}(G,I_0)'_{{\rm
    aff},k}}\chi$ appearing in statement (a) of Theorem \ref{classsusi}. We call the latter ones {\it standard supersingular} ${\mathcal
  H}(G,I_0)_k$-modules. In fact, it is easy to see that these are simple themselves if the pair $(k_{\bullet},{\mathcal
  J})$ satisfies a certain non-periodicity property.

\section{Classical matrix groups}

For $m\in{\mathbb N}$ let $E_m\in{\rm GL}_m$ denote the identity matrix and let ${E^*_d}$ denote the standard antidiagonal element in ${\rm GL}_d$ (i.e. the
permutation matrix of maximal length). Let
$$\widehat{S}_m=\left(\begin{array}{cc}{}&E_m\\-E_m&{}\end{array}\right),\quad\quad\quad {S}_m=\left(\begin{array}{cc}{}&E_m\\E_m&{}\end{array}\right).$$For $1\le i,j\le m$ with $i\ne j$ let $\epsilon_{i,j}\in {\rm GL}_m$ denote the matrix with entries $1$ on the diagonal and at the spot $(i,j)$ and with entries $0$ otherwise.

\subsection{Affine root system $\tilde{C}_d$}

\label{mustersubse}

Assume $d\ge2$. Here $W_{{\rm aff}}$ is the
Coxeter group with generators $s_0, s_1,\ldots,s_d$ and relations
\begin{gather}(s_0s_1)^4=(s_{d-1}s_d)^4=1\quad\quad\quad\mbox{ and }\quad\quad\quad(s_{i-1}s_i)^3=1\quad\mbox{ for }2\le i\le d-1\label{crel}\end{gather}and moreover
$(s_is_j)^2=1$ for all other pairs $i<j$, and $s_i^2=1$ for all $i$. In the extended affine Weyl group $\widehat{W}$
we find (cf. \cite{im}) an element $u$ of length $0$ with \begin{gather}u^2=1\quad\quad\quad\mbox{ and }\quad\quad\quad us_iu=s_{d-i}\quad\mbox{ for }0\le i\le d.\label{ucfo}\end{gather}(We have $\widehat{W}=W_{{\rm aff}}\rtimes W_{\Omega}$ with the two-element subgroup $W_{\Omega}=\{1,u\}$.) 
Consider the symplectic and the general symplectic
group $${\rm Sp}_{2d}({\mathbb Q}_p)=\{A\in {\rm GL}_{2d}({\mathbb
  Q}_p)\,|\,{}^TA\widehat{S}_dA=\widehat{S}_d\},$$$$G={\rm GSp}_{2d}({\mathbb Q}_p)=\{A\in {\rm GL}_{2d}({\mathbb
  Q}_p)\,|\,{}^TA\widehat{S}_dA=\kappa(A)\widehat{S}_d\mbox{ for some }\kappa(A)\in{\mathbb
  Q}_p^{\times}\}.$$Let
$T$ denote the maximal torus in $G$ consisting of all diagonal matrices in
$G$. For $1\le i\le d$ let$$e_i:T\cap {\rm Sp}_{2d}({\mathbb Q}_p)\longrightarrow {\mathbb
   Q}_p^{\times},\quad {\rm
  diag}(x_1,\ldots,x_{2d})\mapsto x_i.$$For $1\le i,j\le d$ and $\epsilon_1, \epsilon_2\in\{\pm 1\}$ we thus obtain characters (using additive notation as usual) $\epsilon_1e_i+\epsilon_2e_j:T\cap {\rm Sp}_{2d}({\mathbb Q}_p)\longrightarrow {\mathbb
   Q}_p^{\times}$. We extend these latter ones to $T$ by setting$$\epsilon_1e_i+\epsilon_2e_j:T\longrightarrow {\mathbb
   Q}_p^{\times},\quad A={\rm
  diag}(x_1,\ldots,x_{2d})\mapsto
x_i^{\epsilon_1}x_j^{\epsilon_2}\kappa(A)^{\frac{-\epsilon_1-\epsilon_2}{2}}.$$For $i=j$ and $\epsilon=\epsilon_1=\epsilon_2$ we simply write $\epsilon 2 e_i$. Then $\Phi=\{\pm
e_i\pm e_j\,|\,i\ne j\}\cup\{\pm 2e_i\}$ is the root system of $G$ with respect to $T$. It is
of type $C_d$. We choose the positive system $\Phi^+=\{
e_i\pm e_j\,|\,i< j\}\cup\{2e_i\,|\,1\le i\le d\}$ with corresponding set of
simple roots $\Delta=\{\alpha_1=e_1-e_2,
\alpha_2=e_2-e_3,\ldots,\alpha_{d-1}=e_{d-1}-e_d, \alpha_{d}=2e_d\}$. The negative of the highest root is $\alpha_0=-2e_1$. For $0\le i\le d$ we have the following explicit formula for $\alpha_i^{\vee}=(\alpha_i)^{\vee}$:\begin{gather}\alpha_i^{\vee}(x)=\left\{\begin{array}{l@{\quad:\quad}l}{\rm diag}(x^{-1},E_{d-1},x,E_{d-1})&\quad
      i=0\\
{\rm diag}(E_{i-1},x,x^{-1},E_{d-i-1},E_{i-1},x^{-1},x,E_{d-i-1})&\quad1\le
i\le d-1\\ 
{\rm diag}(E_{d-1},x,E_{d-1},x^{-1})&\quad i=d\end{array}\right.\label{hsic}\end{gather}

For $\alpha\in \Phi$ let $N_{\alpha}^0$ be the subgroup of the
corresponding root subgroup $N_{\alpha}$ of $G$ all of whose elements belong
to ${\rm GL}_{2d}({\mathbb
  Z}_p)$. Let $I_0$ denote the pro-$p$-Iwahori subgroup generated by the $N_{\alpha}^0$ for
all $\alpha\in \Phi^+$, by the $(N_{\alpha}^0)^p$ for
all $\alpha \in \Phi^-=\Phi-\Phi^+$, and by the maximal pro-$p$-subgroup of $T_0$. Let $I$ denote the Iwahori
subgroup of $G$ containing $I_0$. Let $N_0$ be the subgroup
of $G$ generated by all $N_{\alpha}^0$ for $\alpha\in \Phi^+$.

For $1\le i\le d-1$ define the block diagonal
matrix $$\dot{s}_{i}={\rm
  diag}(E_{i-1},\widehat{S}_1,E_{d-i-1},E_{i-1},\widehat{S}_1,E_{d-i-1})$$and
furthermore$$\dot{s}_d=\left(\begin{array}{cccc}E_{d-1}&{}&{}&{}\\{}&{}&{}&1\\{}&{}&E_{d-1}&{}\\{}&-1&{}&{}\end{array}\right),\quad\quad\quad
\dot{s}_0=\left(\begin{array}{cccc}{}&{}&-p^{-1}&{}\\{}&E_{d-1}&{}&{}\\p&{}&{}&{}\\{}&{}&{}&E_{d-1}\end{array}\right).$$ Then
$\dot{s}_{i}$ for $0\le i\le d$ belongs to ${\rm Sp}_{2d}({\mathbb Q}_p)\subset G$ and normalizes $T$. Its image element $s_{i}=s_{\alpha_i}$ in $N(T)/ZT_0=\widehat{W}$ is the reflection corresponding to $\alpha_i$. The ${s}_0,{s}_1,\ldots,{s}_{d-1},{s}_{d}$ are
Coxeter generators of $W_{\rm aff}\subset \widehat{W}$ satisfying the
relations (\ref{crel}). Put$$\dot{u}=\left(\begin{array}{cccc}{}&{E^*_d}\\p{E^*_d}&{}\end{array}\right).$$Then $\dot{u}$ belongs to $N(T)$ and normalizes $I$ and $I_0$. Its image element $u$ in $N(T)/ZT_0$ satisfies the relations (\ref{ucfo}). In $N(T)$ we consider$$\phi=(p\cdot{\rm id})\dot{s}_{d}\dot{s}_{d-1}\cdots
\dot{s}_{1}\dot{s}_0.$$We may
rewrite this as $\phi=(p\cdot{\rm id})\dot{s}_{\beta(1)}\cdots\dot{s}_{\beta(d+1)}$ where we put
$\beta(i)=d+1-i$ for $1\le i\le d+1$. For $a,b\in{\mathbb Z}_{\ge0}$ with $0\le b<d+1$ we
put$$C^{(a(d+1)+b)}=\phi^as_d\cdots s_{d-b+1}C=\phi^as_{\beta(1)}\cdots s_{\beta(b)}C.$$Define the homomorphism$$\tau:{\mathbb
  Z}_p^{\times}\longrightarrow T_0,\quad\quad x\mapsto {\rm
  diag}(xE_{d},E_{d}).$$

\begin{lem}\label{tauc} We have $\phi^{d}\in T$ and $\phi^{d}N_0\phi^{-d}\subset N_0$. The sequence $C=C^{(0)}, C^{(1)},
C^{(2)},\ldots$ satisfies hypothesis (\ref{n0ineu}). In particular we may define $\alpha^{(j)}\in\Phi^+$
for all $j\ge0$.

(b) For all $j\ge0$ we have $\alpha^{(j)}\circ\tau={\rm id}_{{\mathbb
  Z}_p^{\times}}$. 

(c) We have $\tau(a)\phi=\phi\tau(a)$ for all
$a\in{\mathbb Z}_p^{\times}$.
\end{lem}

{\sc Proof:} (a) This is Lemma \ref{concept} in practice. A matrix computation shows $\phi^{d}=(-1)^{d-1}{\rm
  diag}(p^{d+1}E_{d},p^{d-1}E_{d})\in T$. The group $N_{\alpha}$ for $\alpha\in\Phi^+$ is generated by $\epsilon_{i,j+d}\epsilon_{j,i+d}$ if $\alpha=e_i+e_j$ with $1\le i<j\le d$, by $\epsilon_{i,i+d}$ if $\alpha=2e_i$ with $1\le i\le d$, and by $\epsilon_{i,j}\epsilon^{-1}_{i+d,j+d}$ if $\alpha=e_i-e_j$ with $1\le i<j\le d$. Using this we
find$$\phi^{d}N_0\phi^{-d}=\prod_{\alpha\in\Phi^+}\phi^{d}(N_0\cap
N_{\alpha})\phi^{-d}=\prod_{\alpha\in\Phi^+}(N_0\cap
N_{\alpha})^{p^{m_{\alpha}}},$$$$m_{\alpha}=\left\{\begin{array}{l@{\quad:\quad}l}2&\quad
      \alpha=e_i+e_j\mbox{ with }1\le i<j\le d\\2&\quad\alpha=2e_i\mbox{ with }1\le i\le d\\0&\quad\mbox{all other }\alpha\in\Phi^+\end{array}\right.$$In particular we find $\phi^{d}N_0\phi^{-d}\subset N_0$ and
$[N_0:\phi^{d}N_0\phi^{-d}]=p^{d(d+1)}$. Let $n\in {\mathbb N}$. We have $\ell(\phi^{n})\le{n(d+1)}$ because the image of $\phi$ in $\widehat{W}$ is a product of $d+1$ Coxeter generators. If $m\ge n$ is such that $m\in d{\mathbb N}$ we then have, on the other hand,$$p^{\ell(\phi^{m})}\ge[I_0:(I_0\cap \phi^m I \phi^{-m})]\ge [N_0:(N_0\cap \phi^m I_0 \phi^{-m})]=[N_0:\phi^m N_0 \phi^{-m}]=p^{m(d+1)}.$$Thus $$\ell(\phi^{n})\ge \ell(\phi^{m})-\ell(\phi^{m-n})\ge m(d+1)-(m-n)(d+1)=n(d+1),$$hence $\ell(\phi^{n})={n(d+1)}$. We have shown that $\phi$ is power multiplicative and$$[N_0:\phi^{d}N_0\phi^{-d}]=[N_0:(N_0\cap\phi^{d}N_0\phi^{-d})]=
[I_0:(I_0\cap\phi^{d}I_0\phi^{-d})].$$We get $I_0=N_0\cdot
(I_0\cap\phi^{d}I_0\phi^{-d})$ and
that hypothesis
(\ref{n0ineu}) holds true. 

(b) As $\phi^d\in T$ we have
$\{\alpha^{(j)}\,|\,j\ge0\}=\{\alpha\in\Phi^+\,|\,m_{\alpha}\ne0\}$. This
implies (b). 


(c) Another matrix computation. \hfill$\Box$\\

As explained in subsection \ref{erklaer} we now obtain a functor $M\mapsto {\bf
  D}(\Theta_*{\mathcal V}_M)$ from ${\rm Mod}^{\rm fin}({\mathcal
  H}(G,I_0))$ to the category of
$(\varphi^{d+1},\Gamma)$-modules over ${\mathcal O}_{\mathcal E}$. As explained in \cite{dfun}, to compute
it we need to understand the intermediate objects $H_0(\overline{\mathfrak
  X}_+,\Theta_*{\mathcal V}_M)$, acted on by $\lfloor{\mathfrak
  N}_0,\varphi^{d+1},\Gamma\rfloor$.\\

Let ${\mathcal H}(G,I_0)'_{{\rm
    aff},k}$ denote the $k$-sub algebra of ${\mathcal H}(G,I_0)_{k}$ generated
by ${\mathcal H}(G,I_0)_{{\rm
    aff},k}$ together with $T_{p\cdot {\rm id}}=T_{\dot{u}^{2}}$ and
$T^{-1}_{p\cdot {\rm id}}=T_{p^{-1}\cdot {\rm id}}$.

Suppose we are given a character $\lambda:\overline{T}\to
k^{\times}$, a
subset ${\mathcal J}\subset S_{\lambda}$ and some $b\in k^{\times}$. Define the numbers $0\le k_i=k_i(\lambda,{\mathcal
  J})\le p-1$ as in subsection \ref{basic}. The character $\chi_{\lambda,{\mathcal
  J}}$ of ${\mathcal H}(G,I_0)_{{\rm
    aff},k}$ extends uniquely to a character$$\chi_{\lambda,{\mathcal
  J},b}:{\mathcal H}(G,I_0)'_{{\rm
    aff},k}\longrightarrow k$$which sends $T_{p\cdot {\rm id}}$ to $b$ (see the
proof of \cite{vigneras} Proposition 3). Define the ${\mathcal H}(G,I_0)_k$-module$$M=M[\lambda,{\mathcal
  J},b]={\mathcal H}(G,I_0)_k\otimes_{{\mathcal H}(G,I_0)'_{{\rm
    aff},k}}k.e$$where $k.e$ denotes the one dimensional $k$-vector space on
the basis element $e$, endowed with the action of ${\mathcal H}(G,I_0)'_{{\rm
    aff},k}$ by the character $\chi_{\lambda,{\mathcal
  J},b}$. As a $k$-vector space, $M$ has dimension $2$, a $k$-basis is
$e,f$ where we write $e=1\otimes e$ and $f=T_{\dot{u}}\otimes e$.\\

{\bf Definition:} We call an ${\mathcal
  H}(G,I_0)_k$-module {\it standard supersingular} if it is isomorphic with $M[\lambda,{\mathcal
  J},b]$ for some $\lambda, {\mathcal
  J}, b$ such that $k_{\bullet}\notin\{(0,\ldots,0), (p-1,\ldots,p-1)\}$.\\

For $0\le j\le d$ put $\widetilde{j}=d-j$. Letting
$\widetilde{\beta}=\widetilde{(.)}\circ\beta$ we then have
$$\dot{u}\phi\dot{u}^{-1}=(p\cdot{\rm id})\dot{s}_{\widetilde{\beta}(1)}\cdots\dot{s}_{\widetilde{\beta}(d+1)}.$$Put $n_e=\sum_{i=0}^{d}k_{d-i}p^i=\sum_{i=0}^{d}k_{\beta(i+1)}p^i$ and $n_f=\sum_{i=0}^{d}k_{i}p^i=\sum_{i=0}^{d}k_{\widetilde{\beta}(i+1)}p^i$. Put
$\varrho=\prod_{i=0}^{d}(k_{i}!)=\prod_{i=0}^d(k_{\beta(i+1)}!)=\prod_{i=0}^d(k_{\widetilde{\beta}(i+1)}!)$. Let $0\le s_e, s_f\le p-2$ be such that
$\lambda(\tau(x))=x^{-s_e}$ and $\lambda(\dot{u}\tau(x)\dot{u}^{-1})=x^{-s_f}$
for all $x\in{\mathbb F}_p^{\times}$.

\begin{lem}\label{combisusic} The assignment $M[\lambda,{\mathcal
  J},b]\mapsto (n_e,s_e,b\varrho^{-1})$ induces a bijection
between the set of isomorphism classes of standard supersingular ${\mathcal
  H}(G,I_0)_k$-modules and ${\mathfrak S}_C(d+1)$. 
\end{lem}

{\sc Proof:} We have $\prod_{i=0}^d\alpha^{\vee}_{i}(x)=1$ for all
$x\in{\mathbb F}_p^{\times}$ (as can be seen e.g. from formula
(\ref{hsic})). This implies \begin{gather}\sum_{i=0}^dk_i\equiv n_e\equiv n_f\equiv 0\quad\mbox{ mod }(p-1).\label{perrelc}\end{gather}One can deduce from \cite{vigneras} Proposition 3 that for two sets of data $\lambda, {\mathcal
  J}, b$ and $\lambda', {\mathcal
  J}', b'$ the ${\mathcal H}(G,I_0)_k$-modules $M[\lambda,{\mathcal
  J},b]$ and $M[\lambda',{\mathcal
  J}',b']$ are isomorphic if and only if $b=b'$ and the pair $(\lambda, {\mathcal
  J})$ is conjugate with the pair $(\lambda', {\mathcal
  J}')$ by means of a power of $\dot{u}$, i.e. by means of $\dot{u}^0=1$ or
$\dot{u}^1=\dot{u}$. Conjugating
$(\lambda,{\mathcal J})$ by $\dot{u}$ has the effect of substituting
$k_{d-i}$ with $k_i$, for any $i$. The datum of the character $\lambda$ is
equivalent with the datum of $s_e$ together with all the $k_i$ taken modulo $(p-1)$ since the images of $\tau$ and all
$\alpha_i^{\vee}$ together generate $\overline{T}$. Knowing the set ${\mathcal
  J}$ is then equivalent with knowing the numbers $k_i$ themselves (not just modulo $(p-1)$). Thus, our mapping is
well defined and bijective. \hfill$\Box$\\

Fix $M=M[\lambda,{\mathcal
  J},b]$. Let $0\le j\le d$. Normalize the homomorphism $\iota_{\alpha_j}:{\rm
  SL}_2({\mathbb Q}_p)\to G$ such that $\iota_{\alpha_j}(\dot{s})=\dot{s}_j$
for $\dot{s}=\left( \begin{array}{cc} 0 & -1 \\
    1 & 0\end{array}\right)$. Let $t_j=\iota_{\alpha_j}([\nu])-1\in k[[\iota_{\alpha_j}{\mathfrak N}_0]]\subset k[[N_0]]$. Let $F_j$ denote the codimension-$1$-face of $C$ contained in the (affine) reflection hyperplane (in $A\subset X$) for $s_j$. 

\begin{lem}\label{locca} In ${\mathcal V}_M^X(F_j)$ we have
  $t_j^{k_j}\dot{s}_je=k_j!e$ and $t_j^{k_{d-j}}\dot{s}_jf=k_{d-j}!f$ for all $0\le j\le d$.
\end{lem}

{\sc Proof:} Let $I_0^{j}$ denote the subgroup of $G$ generated
by $I_0$ and $\dot{s}_{j}I_0 \dot{s}_{j}^{-1}$. Let $\overline{I}_0^{j}$ denote the maximal reductive (over ${\mathbb F}_p$)
quotient of $I_0^{j}$. Then $\iota_{\alpha_j}({\rm SL}_2({\mathbb
  Z}_p))\subset I_0^{j}$ and $\iota_{\alpha_j}$ induces an isomorphism between
either the quotient $\overline{\mathcal S}={\rm SL}_2({\mathbb F}_p)$ of ${\rm SL}_2({\mathbb
  Z}_p)$, or the quotient $\overline{\mathcal S}={\rm PSL}_2({\mathbb F}_p)$ of ${\rm SL}_2({\mathbb
  Z}_p)$, with the quotient $\overline{I}_0^{j}$ of $I_0^{j}$. The
subgroup $\overline{\mathcal U}$ of $\overline{\mathcal S}$ generated by (the image of) $\nu$
maps isomorphically to the image of $I_0$ in
$\overline{I}_0^{j}$. Therefore we get an embedding of $k$-algebras\begin{align}{\mathcal
  H}(\overline{\mathcal S},\overline{\mathcal U})_k={\rm End}_{k[\overline{\mathcal S}]}({\rm
  ind}_{\overline{\mathcal U}}^{\overline{\mathcal S}}{\bf 1}_k)^{\rm
  op}&\cong{\rm End}_{k[I_0^{j}]}({\rm
  ind}_{I_0}^{I_0^{j}}{\bf 1}_k)^{\rm op}\notag\\{}&\hookrightarrow{\mathcal H}(G,I_0)_{{\rm
    aff},{k}}\quad\subset\quad {\mathcal
  H}(G,I_0)_{k}.\notag\end{align}Let
$\rho:{\mathcal
  H}(\overline{\mathcal S},\overline{\mathcal U})_k\to k$ be the character
defined by $\rho(T_{\dot{s}})=-1$ if $k_{j}=p-1$ but
$\rho(T_{\dot{s}})=0$ if $0\le k_{j}<p-1$, and by
$\rho(T_{h(x)})=x^{-k_{j}}$. We have a morphism of ${\rm
  SL}_2({\mathbb Z}_p)$-representations (acting on ${\mathcal V}_M^X(F_j)$ through $\iota_{\alpha_j}$)$$({\rm ind}_{\overline{\mathcal
    U}}^{\overline{\mathcal S}}{\bf 1}_k)\otimes_{{\mathcal
  H}(\overline{\mathcal S},\overline{\mathcal U})_k}\rho\longrightarrow
{\mathcal V}_M^X(F_j)$$sending a basis element $e$ of $\rho$ to $e\in{\mathcal V}_M^X(F_j)$. It
is injective since its restriction to the space $k.e$ of ${\mathcal
  U}$-invariants (see \cite{dfun} Lemma 2.3) is injective. By \cite{dfun}
Lemma 2.5 we have
$t^{k_{j}}\dot{s}e=k_{j}!e$ in $({\rm ind}_{\overline{\mathcal
    U}}^{\overline{\mathcal S}}{\bf 1}_k)\otimes_{{\mathcal
  H}(\overline{\mathcal S},\overline{\mathcal U})_k}\rho$ where $t=[\nu]-1\in
k[\overline{\mathcal
    U}]$. We get $t_j^{k_j}\dot{s}_je=k_j!e$ in ${\mathcal V}_M^X(F_j)$. A
  similar argument gives $t_j^{k_{d-j}}\dot{s}_jf=k_{d-j}!f$. For this notice that as the Hecke operator $T_t$ for
$t\in\overline{T}$ acts on $k.e$ through $\lambda(t^{-1})$, it acts on $k.f=k.T_{\dot{u}}e$
through $\lambda (\dot{u}t^{-1}\dot{u}^{-1})$ (the same computation as in
formula (\ref{exempel}) below), and that formula (\ref{kljfmlc}) implies $\lambda (\dot{u}\alpha_j^{\vee}(x)\dot{u}^{-1})=x^{k_{d-j}}$.\hfill$\Box$\\

\begin{lem}\label{globca} For an appropriate choice of the isomorphism
  $\Theta$ we have, in $H_0(\overline{\mathfrak X}_+,\Theta_*{\mathcal V}_M)$,
\begin{gather}t^{n_e}\varphi^{d+1}e=\varrho
  b^{-1}e,\label{cstaef}\\t^{n_f}\varphi^{d+1} f=\varrho
  b^{-1}f,\label{cstafe}\\\gamma(x)e=x^{-s_e}e,\label{gaac1}\\\gamma(x)f=x^{-s_f}f\label{gaac2}\end{gather}for
$x\in{\mathbb F}_p^{\times}$. The action of $\Gamma_0$ on
$H_0(\overline{\mathfrak X}_+,\Theta_*{\mathcal V}_M)$ is trivial on the subspace $M$.
\end{lem} 

{\sc Proof:} The construction of $\Theta$ in \cite{dfun} is based on the choice
 of elements $\nu_i\in N_0$ topologically generating the relevant pro-$p$ subgroup 'at
 distance $i$', for each $0\le i\le d$. These $\nu_i$ can be taken to be the elements
 $y_{i-1}\cdots y_0\iota_{d-i}(\nu)y_0^{-1}\cdots y_{i-1}^{-1}$. Their actions on the halftree $Y$ correspond to the actions
 of the $\nu^{p^{i}}$ on the halftree ${\mathfrak X}_+$. 

We use the notations and the statements of Lemma
\ref{prepa}, observing $\beta(i+1)=d-i$. For $0\le i\le d$ we have $y_{i-1}\cdots
y_0F_{d-i}=C^{(i)}\cap C^{(i+1)}$ and $y_{i-1}\cdots y_0C=C^{(i)}$. Thus $y_{i-1}\cdots y_0$ defines isomorphisms $${\mathcal V}_M^X(F_{d-i})\cong
 {\mathcal V}^X_M(C^{(i)}\cap C^{(i+1)})\quad\quad\mbox{ and }\quad\quad {\mathcal V}_M^X(C)\cong
 {\mathcal V}^X_M(C^{(i)}).$$But we have $${\mathcal V}^X_M(C^{(i)}\cap C^{(i+1)})=\Theta_*{\mathcal V}_M({\mathfrak
   v}_i)\quad\quad\mbox{ and }\quad\quad{\mathcal V}^X_M(C^{(i)})=\Theta_*{\mathcal V}_M({\mathfrak
   e}_i).$$It follows that, under the above isomorphisms, the action of $t_{d-i}$, resp. of
 $\dot{s}_{d-i}$, on ${\mathcal V}_M^X(F_{d-i})$ becomes the action of
 $[\nu]^{p^{i}}-1$, resp. of $y_i$, on $\Theta_*{\mathcal V}_M({\mathfrak
   v}_i)$. Now as we are in characteristic $p$ we have
 $t^{p^i}=([\nu]-1)^{p^{i}}=[\nu]^{p^{i}}-1$. Applying this to the element
 $e$, resp. $f$, of ${\mathcal
  V}_M^X(C)\subset{\mathcal V}_M^X(F_{d-i})$, Lemma \ref{locca}
 tells us $$(t^{p^i})^{k_{d-i}}y_i\cdots y_0e=k_{d-i}! y_{i-1}\cdots
 y_0e\quad\mbox{ resp. }\quad(t^{p^i})^{k_{i}}y_i\cdots y_0f=k_{i}! y_{i-1}\cdots
 y_0f.$$We compose these formulae for all $0\le i\le d$ and finally recall that the
 central element $p\cdot{\rm id}$ acts on $M$ through the Hecke operator
 $T_{p^{-1}\cdot{\rm id}}$, i.e. by $b^{-1}$. We get formulae (\ref{cstafe}) and (\ref{cstaef}). 

Next recall that the action of
$\gamma(x)$ on $H_0(\overline{\mathfrak X}_+,\Theta_*{\mathcal V}_M)$ is given
by that of $\tau(x)\in T$, i.e. by the Hecke operator $T_{\tau(x)^{-1}}$. We thus
compute \begin{gather}\gamma(x)e=T_{\tau(x)^{-1}}e=\lambda(\tau(x))e,\notag\\\gamma(x)f=T_{\tau(x)^{-1}}T_{\dot{u}}e=T_{\dot{u}\tau(x)^{-1}}e=T_{\dot{u}}
  T_{\dot{u}\tau(x)^{-1}\dot{u}^{-1}}e=T_{\dot{u}}\lambda(\dot{u}\tau(x)\dot{u}^{-1})e=\lambda(\dot{u}\tau(x)\dot{u}^{-1})f\label{exempel}\end{gather}and
obtain formulae (\ref{gaac1}) and (\ref{gaac2}).\hfill$\Box$\\

\begin{kor}\label{cphigaa} The \'{e}tale $(\varphi^{d+1},\Gamma)$-module ${\bf
  D}(\Theta_*{\mathcal V}_M)$ over $k_{\mathcal E}$ associated with $H_0(\overline{\mathfrak X}_+,\Theta_*{\mathcal V}_M)$ admits a $k_{\mathcal E}$-basis $g_e$, $g_f$ such that \begin{gather}\varphi^{d+1}g_e=b\varrho^{-1}t^{n_e+1-p^{d+1}}g_e\notag\\\varphi^{d+1}g_f=b\varrho^{-1}t^{n_f+1-p^{d+1}}g_f\notag\\\gamma(x)g_e-x^{s_e}g_e\in t\cdot k_{\mathcal E}^+\cdot g_e\notag\\\gamma(x)g_f-x^{s_f}g_f\in t\cdot k_{\mathcal E}^+\cdot g_f.\notag\end{gather}
\end{kor}

{\sc Proof:} This follows from Lemma \ref{globca} as explained in \cite{dfun}
Lemma 6.4.\hfill$\Box$\\ 

\begin{kor}\label{cflip} The functor $M\mapsto {\bf
  D}(\Theta_*{\mathcal V}_M)$ induces a bijection between

(a) the set of isomorphism classes of standard supersingular ${\mathcal
  H}(G,I_0)_k$-modules, and

(b) the set of isomorphism classes of $C$-symmetric \'{e}tale
$(\varphi^{d+1},\Gamma)$-modules over $k_{\mathcal E}$.  
\end{kor}

{\sc Proof:} For $x\in{\mathbb F}_p^{\times}$ we have $$\tau(x)\cdot \dot{u}\tau^{-1}(x)\dot{u}^{-1}={\rm
  diag}(xE_d,x^{-1}E_d)=(\sum_{i=0}^d(i+1)\alpha^{\vee}_i)(x)$$in $\overline{T}$. Applying
$\lambda$ and observing $\sum_{i=0}^dk_i\equiv 0$ modulo $(p-1)$ we get
$$x^{s_f-s_e}=\lambda((\sum_{i=0}^d(i+1)\alpha^{\vee}_i)(x))=x^{\sum_{i=0}^dik_i}$$and
hence $s_f-s_e\equiv \sum_{i=0}^dik_i$ modulo $(p-1)$. Therefore $${\bf
  D}(\Theta_*{\mathcal V}_M)=  {\bf
  D}(\Theta_*{\mathcal V}_M)_1\oplus {\bf
  D}(\Theta_*{\mathcal V}_M)_2$$with ${\bf
  D}(\Theta_*{\mathcal V}_M)_1=k_{\mathcal E}g_e$ and ${\bf
  D}(\Theta_*{\mathcal V}_M)_2=k_{\mathcal E}g_f$ is a direct sum decomposition of ${\bf
  D}(\Theta_*{\mathcal V}_M)$ as an \'{e}tale
$(\varphi^{d+1},\Gamma)$-module over $k_{\mathcal E}$, and by Corollary (\ref{cphigaa}) this decomposition identifies ${\bf
  D}(\Theta_*{\mathcal V}_M)$ as being $C$-symmetric (recall in particular that $n_e=\sum_{i=0}^{d}k_{d-i}p^i$ and $n_f=\sum_{i=0}^{d}k_{i}p^i$). Composing this assignment $M\mapsto {\bf
  D}(\Theta_*{\mathcal V}_M)$ (from standard supersingular ${\mathcal
  H}(G,I_0)_k$-modules to $C$-symmetric \'{e}tale
$(\varphi^{d+1},\Gamma)$-modules) with the bijection of Lemma \ref{combiflip} (i) (from $C$-symmetric \'{e}tale
$(\varphi^{d+1},\Gamma)$-modules to ${\mathfrak S}_C(d+1)$) we get the bijection of Lemma \ref{combisusic} (from standard supersingular ${\mathcal
  H}(G,I_0)_k$-modules to ${\mathfrak S}_C(d+1)$). Thus, $M\mapsto {\bf
  D}(\Theta_*{\mathcal V}_M)$ is a bijection as stated.\hfill$\Box$\\ 

{\bf Remark:} Consider the subgroup $G'={\rm Sp}_{2d}({\mathbb Q}_p)$ of $G$. If we replace the above $\tau$ by $\tau:{\mathbb
  Z}_p^{\times}\longrightarrow T_0$, $x\mapsto {\rm
  diag}(xE_{d},x^{-1}E_{d})$ and if we replace the above $\phi$ by $\phi=\dot{s}_{d}\dot{s}_{d-1}\cdots
\dot{s}_{1}\dot{s}_0$ then everything in fact happens inside
$G'$. We then have $\alpha^{(j)}\circ\tau={\rm id}^2_{{\mathbb
  Z}_p^{\times}}$ for all $j\ge0$. Let ${\rm Mod}_0^{\rm fin}{\mathcal
  H}(G',G'\cap I_0)$ denote the category of
finite-${\mathfrak o}$-length ${\mathcal
  H}(G',G'\cap I_0)$-modules on which $\tau(-1)$ (i.e. $T_{\tau(-1)}=T_{\tau(-1)^{-1}}$) acts trivially. For $M\in{\rm Mod}_0^{\rm fin}{\mathcal
  H}(G',G'\cap I_0)$ we obtain an action of $\lfloor{\mathfrak N}_0,\varphi^{d+1},\Gamma^2\rfloor$ on
$H_0(\overline{\mathfrak X}_+,\Theta_*{\mathcal V}_M)$, where
$\Gamma^2=\{\gamma^2\,|\,\gamma\in\Gamma\}\subset\Gamma$. Correspondingly,
following \cite{dfun} (as a slight variation from what we explained in
subsection \ref{erklaer}), we obtain a functor from ${\rm Mod}^{\rm fin}_0{\mathcal
  H}(G',G'\cap I_0)$ to the
category of $(\varphi^{d+1},\Gamma^2)$-modules over ${\mathcal O}_{\mathcal
  E}$.\\

{\bf Remark:} In the case $d=2$ one may also work with 
$\phi=(p\cdot{\rm id})\dot{s}_2\dot{s}_1\dot{s}_2\dot{u}$. Its square is
$p\cdot{\rm id}$ times the
square of the $\phi=(p\cdot{\rm id})\dot{s}_2\dot{s}_1\dot{s}_0$ used above.\\

\subsection{Affine root system $\tilde{B}_d$} 

\label{subsecb}

Assume $d\ge 3$. Here $W_{{\rm aff}}$ is the
Coxeter group with generators $s_0, s_1,\ldots,s_d$ and relations
\begin{gather}(s_{d-1}s_{d})^4=1\quad\quad\quad\mbox{ and }\quad\quad\quad(s_{0}s_2)^3=(s_{i-1}s_i)^3=1\quad\mbox{ for }2\le i\le d-1\label{brel}\end{gather}and moreover
$(s_is_j)^2=1$ for all other pairs $i<j$, and $s_i^2=1$ for all $i$. In the extended affine Weyl group $\widehat{W}$
we find (cf. \cite{im}) an element $u$ of length $0$ with \begin{gather}u^2=1\quad\quad\quad\mbox{ and
}\quad\quad\quad us_0u=s_{1}\quad\quad\quad\mbox{ and
}\quad\quad\quad     us_iu=s_{i}\quad\mbox{ for }2\le i\le d.\label{ubfo}\end{gather}
(We have $\widehat{W}=W_{{\rm aff}}\rtimes W_{\Omega}$ with the two-element subgroup $W_{\Omega}=\{1,u\}$.) 
Let$$\widetilde{S}_d=\left(\begin{array}{cc}S_d&0\\0&1\end{array}\right)\in {\rm GL}_{2d+1}({\mathbb
  Q}_p)$$and consider the special orthogonal group $$G={\rm SO}_{2d+1}({\mathbb Q}_p)=\{A\in {\rm SL}_{2d+1}({\mathbb
  Q}_p)\,|\,{}^TA\widetilde{S}_dA=\widetilde{S}_d\}.$$Let
$T$ denote the maximal torus consisting of all diagonal matrices in
$G$. For $1\le i\le d$ let $$e_i:T\longrightarrow {\mathbb Q}_p^{\times},\quad{\rm
  diag}(x_1,\ldots,x_d,x_1^{-1},\ldots,x_d^{-1},1)\mapsto x_i.$$Then (in {\it additive} notation) $\Phi=\{\pm
e_i\pm e_j\,|\,i\ne j\}\cup\{\pm e_i\}$ is the root system of $G$ with respect to $T$. It is
of type $B_d$. We choose the positive system $\Phi^+=\{
e_i\pm e_j\,|\,i< j\}\cup\{e_i\,|\,1\le i\le d\}$ with corresponding set of
simple roots $\Delta=\{\alpha_1=e_1-e_2,
\alpha_2=e_2-e_3,\ldots,\alpha_{d-1}=e_{d-1}-e_d, \alpha_{d}=e_d\}$. The negative of the highest root is $\alpha_0=-e_1-e_2$. For $0\le i\le d$ we have the following explicit formula for $\alpha_i^{\vee}=(\alpha_i)^{\vee}$:\begin{gather}\alpha_i^{\vee}(x)=\left\{\begin{array}{l@{\quad:\quad}l}{\rm diag}(x^{-1},x^{-1},E_{d-2},x,x,E_{d-2},1)&\quad
      i=0\\
{\rm
  diag}(E_{i-1},x,x^{-1},E_{d-i-1},E_{i-1},x^{-1},x,E_{d-i-1},1)&\quad 1\le
i\le d-1\\ 
{\rm diag}(E_{d-1},x^2,E_{d-1},x^{-2},1)&\quad
i=d\end{array}\right.\label{hsib}\end{gather}

{\bf Remark:} For roots $\alpha\in\Phi$ of the form $\alpha=\pm e_i\pm e_j$ the homomorphism
$\iota_{\alpha}:{\rm SL}_2\to {\rm SO}_{2d+1}$ is injective. For roots
$\alpha\in\Phi$ of the form $\alpha=\pm e_i$ the homomorphism
$\iota_{\alpha}:{\rm SL}_2\to {\rm SO}_{2d+1}$ induces an embedding ${\rm
  PSL}_2\to {\rm SO}_{2d+1}$.\\

 For $\alpha\in \Phi$ let $N_{\alpha}^0$ be the subgroup of the
corresponding root subgroup $N_{\alpha}$ of $G$ all of whose elements belong
to ${\rm SL}_{2d+1}({\mathbb
  Z}_p)$. Let $I_0$ denote the pro-$p$-Iwahori subgroup generated by the $N_{\alpha}^0$ for
all $\alpha\in \Phi^+$, by the $(N_{\alpha}^0)^p$ for
all $\alpha \in \Phi^-=\Phi-\Phi^+$, and by the maximal pro-$p$-subgroup of $T_0$. Let $I$ denote the Iwahori
subgroup of $G$ containing $I_0$. Let $N_0$ be the subgroup
of $G$ generated by all $N_{\alpha}^0$ for $\alpha\in \Phi^+$.

For $1\le i\le d-1$ define the block diagonal
matrix $$\dot{s}_i={\rm
  diag}(E_{i-1},{S}_1,E_{d-i-1},E_{i-1},{S}_1,E_{d-i-1},1)$$and furthermore$$\dot{s}_d=\left(\begin{array}{ccccc}E_{d-1}&{}&{}&{}&{}\\{}&{}&{}&1&{}\\{}&{}&E_{d-1}&{}&{}\\{}&1&{}&{}&{}\\{}&{}&{}&{}&-1\end{array}\right).$$Define$$\dot{u}=\left(\begin{array}{ccccc}{}&{}&p^{-1}&{}&{}\\{}&E_{d-1}&{}&{}&{}\\p&{}&{}&{}&{}\\{}&{}&{}&E_{d-1}&{}\\{}&{}&{}&{}&-1\end{array}\right)$$and
$\dot{s}_0=\dot{u}\dot{s}_1\dot{u}$. Then $\dot{s}_{i}$ for $0\le i\le d$ belongs to $G$ and normalizes $T$. Its image element $s_{i}=s_{\alpha_i}$ in $N(T)/T_0=\widehat{W}$ is the reflection corresponding to $\alpha_i$. The ${s}_0,{s}_1,\ldots,{s}_{d-1},{s}_{d}$ are
Coxeter generators of $W_{\rm aff}\subset \widehat{W}$ satisfying the relations (\ref{brel}). Also $\dot{u}$ belongs to $G$; it normalizes $T$, $I$ and $I_0$. Its image element $u$ in $N(T)/T_0=\widehat{W}$ satisfies the relations (\ref{ubfo}). In $N(T)$ we
consider the
element\begin{gather}\phi=\dot{s}_1\dot{s}_2\cdots\dot{s}_{d-1}\dot{s}_{d}\dot{s}_{d-1}\cdots\dot{s}_2\dot{s}_0.\label{phicab}\end{gather}We may
rewrite this as $\phi=\dot{s}_{\beta(1)}\cdots\dot{s}_{\beta(2d-1)}$ where we put
$\beta(i)=i$ for $1\le i\le d$ and $\beta(i)=2d-i$ for $d\le i\le 2d-2$ and
$\beta(2d-1)=0$. We
put$$C^{(a(2d-1)+b)}=\phi^as_{\beta(1)}\cdots s_{\beta(b)}C$$for $a,b\in{\mathbb Z}_{\ge0}$ with $0\le b<2d-1$. Define the homomorphism$$\tau:{\mathbb
  Z}_p^{\times}\longrightarrow T_0,\quad\quad x\mapsto {\rm
  diag}(x,E_{d-1},x^{-1},E_{d-1},1).$$

\begin{lem}\label{taub} We have $\phi^{2}\in T$ and $\phi^{2}N_0\phi^{-2}\subset N_0$. The sequence $C=C^{(0)}, C^{(1)},
C^{(2)},\ldots$ satisfies hypothesis (\ref{n0ineu}). In particular we may define $\alpha^{(j)}\in\Phi^+$
for all $j\ge0$. 

(b) For any $j\ge0$ we have $\alpha^{(j)}\circ\tau={\rm id}_{{\mathbb
  Z}_p^{\times}}$. 

(c) We have $\tau(a)\phi=\phi\tau(a)$ for all
$a\in{\mathbb Z}_p^{\times}$.
\end{lem}

{\sc Proof:} (a) A matrix computation shows $\phi^{2}={\rm
  diag}(p^2,E_{d-1},p^{-2},E_{d-1},1)$. The group $N_{\alpha}$ for $\alpha\in\Phi^+$ is generated by $\epsilon_{i,j+d}\epsilon^{-1}_{j,i+d}$ if $\alpha=e_i+e_j$ with $1\le i<j\le d$, by $\epsilon_{i,j}\epsilon^{-1}_{j+d,i+d}$ if $\alpha=e_i-e_j$ with $1\le i<j\le d$, and by $\epsilon_{i,2d+1}\epsilon^{-1}_{2d+1,i+d}$ if $\alpha=e_i$ with $1\le i\le d$. Using this we
find$$\phi^{2}N_0\phi^{-2}=\prod_{\alpha\in\Phi^+}\phi^{2}(N_0\cap
N_{\alpha})\phi^{-2}=\prod_{\alpha\in\Phi^+}(N_0\cap
N_{\alpha})^{p^{m_{\alpha}}},$$$$m_{\alpha}=\left\{\begin{array}{l@{\quad:\quad}l}2&\quad
      \alpha=e_1-e_i\mbox{ with }1< i\\
2&\quad\alpha=e_1+e_i\mbox{ with }1<i\\2&\quad\alpha=e_1\\
0&\quad\mbox{all other }\alpha\in\Phi^+\end{array}\right.$$In particular we find $\phi^{2}N_0\phi^{-2}\subset N_0$ and
$[N_0:\phi^{2}N_0\phi^{-2}]=p^{2(2d-1)}$. On the other hand, the image of $\phi$ in $\widehat{W}$ is a product of $2d-1$ Coxeter generators. Arguing as in the proof of Lemma \ref{tauc} we combine these facts to obtain $\ell(\phi)=2d-1$ and that $\phi$ is power multiplicative. We obtain that hypothesis
(\ref{n0ineu}) holds true, again by the same reasoning as in Lemma \ref{tauc}. 

(b) As $\phi^{2}\in T$ we have
$\{\alpha^{(j)}\,|\,j\ge0\}=\{\alpha\in\Phi^+\,|\,m_{\alpha}\ne0\}$. This
implies (b). 


(c) Another matrix computation.\hfill$\Box$\\

As explained in subsection \ref{erklaer} we now obtain a functor from ${\rm Mod}^{\rm fin}({\mathcal
  H}(G,I_0))$ to the category of
$(\varphi^{2d-1},\Gamma)$-modules over ${\mathcal O}_{\mathcal E}$.\\

Suppose we are given a character $\lambda:\overline{T}\to
k^{\times}$ and a
subset ${\mathcal J}\subset S_{\lambda}$. Define the numbers $0\le k_i=k_i(\lambda,{\mathcal
  J})\le p-1$ as in subsection \ref{basic}. Notice that $k_d$ is necessarily even since in $\alpha_d^{\vee}(x)$ (formula (\ref{hsib})) the entry $x$ only appears squared. Define the ${\mathcal H}(G,I_0)_k$-module$$M=M[\lambda,{\mathcal
  J}]={\mathcal H}(G,I_0)_k\otimes_{{\mathcal H}(G,I_0)_{{\rm
    aff},k}}k.e$$where $k.e$ denotes the one dimensional $k$-vector space on
the basis element $e$, endowed with the action of ${\mathcal H}(G,I_0)_{{\rm
    aff},k}$ by the character $\chi_{\lambda,{\mathcal
  J}}$. As a $k$-vector space, $M$ has dimension $2$, a $k$-basis is
$e,f$ where we write $e=1\otimes e$ and $f=T_{\dot{u}}\otimes e$.\\

{\bf Definition:} We call an ${\mathcal
  H}(G,I_0)_k$-module {\it standard supersingular} if it is isomorphic with $M[\lambda,{\mathcal
  J}]$ for some $\lambda, {\mathcal
  J}$ such that $k_{\bullet}\notin\{(0,\ldots,0), (p-1,\ldots,p-1)\}$.

 A {\it packet} of standard supersingular ${\mathcal
  H}(G,I_0)_k$-modules is a set of isomorphism classes of standard supersingular ${\mathcal
  H}(G,I_0)_k$-modules all of which give rise to the same ${\mathcal
  J}$ and the same $k_{\bullet}=(k_i)_{i}$.  \\

{\bf Remark:} Let $\overline{T}'$ denote the subgroup of $\overline{T}$ generated by the $\alpha_i^{\vee}({\mathbb F}_p^{\times})$ for all $0\le i\le d$. Formula (\ref{hsib}) implies $[\overline{T}:\overline{T}']=2$. Two supersingular ${\mathcal
  H}(G,I_0)_k$-modules $M[\lambda,{\mathcal
  J}]$ and $M[\lambda',{\mathcal
  J}']$ belong to the same packet if and only if ${\mathcal
  J}={\mathcal
  J}'$ and if the restrictions of $\lambda$ and $\lambda'$ to $\overline{T}'$ coincide.\\

For $2\le j\le d$ we put $\widetilde{j}=j$, furthermore we put
$\widetilde{0}=1$ and $\widetilde{1}=0$. Letting
$\widetilde{\beta}=\widetilde{(.)}\circ\beta$ we then have
$$\dot{u}\phi\dot{u}^{-1}=\dot{s}_{\widetilde{\beta}(1)}\cdots\dot{s}_{\widetilde{\beta}(2d-1)}.$$Put $n_e=\sum_{i=0}^{2d-2}k_{\beta(i+1)}p^i$ and
$n_f=\sum_{i=0}^{2d-2}k_{\widetilde{\beta}(i+1)}p^i$. Put
$\varrho=k_0!k_1!k_d!\prod_{i=2}^{d-1}(k_{i}!)^2=\prod_{i=0}^{2d-2}(k_{\beta(i+1)}!)=\prod_{i=0}^{2d-2}(k_{\widetilde{\beta}(i+1)}!)$. Let $0\le s_e, s_f\le p-2$
be such that
$\lambda(\tau(x))=x^{-s_e}$ and $\lambda(\dot{u}\tau(x)\dot{u}^{-1})=x^{-s_f}$
for all $x\in{\mathbb F}_p^{\times}$.

\begin{lem}\label{combisusib} The assignment $M[\lambda,{\mathcal
  J}]\mapsto (n_e,s_e)$ induces a bijection
between the set of packets of standard supersingular ${\mathcal
  H}(G,I_0)_k$-modules and ${\mathfrak S}_B(2d-1)$. 
\end{lem}

{\sc Proof:} We have $\alpha^{\vee}_{0}(x)\alpha^{\vee}_{1}(x)\alpha^{\vee}_{{d}}(x)\prod_{i=2}^{d-1}(\alpha_{i}^{\vee})^2(x)=1$ for all $x\in{\mathbb F}_p^{\times}$ (as can be seen e.g. from formula (\ref{hsib})). This implies\begin{gather}k_0+k_1+k_d+2\sum_{i=2}^{d-1}k_i\equiv
  n_e\equiv n_f\equiv 0\quad\mbox{ mod }(p-1).\label{perrelb}\end{gather}We
further proceed exactly as in the proof of Lemma \ref{combisusic}.\hfill$\Box$\\

\begin{lem}\label{globba} Let $M=M[\lambda,{\mathcal
  J}]$ for some $\lambda$, ${\mathcal
  J}$. For an appropriate choice of the isomorphism
  $\Theta$ we have, in $H_0(\overline{\mathfrak X}_+,\Theta_*{\mathcal V}_M)$,
\begin{gather}t^{n_e}\varphi^{2d-1}e=\varrho
  e,\label{bstaef}\\t^{n_f}\varphi^{2d-1} f=\varrho
  f,\label{bstafe}\\\gamma(x)e=x^{-s_e}e,\label{gaab1}\\\gamma(x)f=x^{-s_f}f\label{gaab2}\end{gather}for
$x\in{\mathbb F}_p^{\times}$. The action of $\Gamma_0$ on
$H_0(\overline{\mathfrak X}_+,\Theta_*{\mathcal V}_M)$ is trivial on the subspace $M$.
\end{lem} 

{\sc Proof:} As in Lemma \ref{globca}. \hfill$\Box$\\

\begin{kor}\label{cphiga} The \'{e}tale $(\varphi^{2d-1},\Gamma)$-module over $k_{\mathcal E}$ associated with $H_0(\overline{\mathfrak X}_+,\Theta_*{\mathcal V}_M)$ admits a $k_{\mathcal E}$-basis $g_e$, $g_f$ such that \begin{gather}\varphi^{2d-1}g_e=\varrho^{-1}t^{n_e+1-p^{2d-1}}g_e\notag\\\varphi^{2d-1}g_f=\varrho^{-1}t^{n_f+1-p^{2d-1}}g_f\notag\\\gamma(x)g_e-x^{s_e}g_e\in t\cdot k_{\mathcal E}^+\cdot g_e\notag\\\gamma(x)g_f-x^{s_f}g_f\in t\cdot k_{\mathcal E}^+\cdot g_f.\notag\end{gather}
\end{kor}

{\sc Proof:} This follows from Lemma \ref{globba} as explained in \cite{dfun}
Lemma 6.4.\hfill$\Box$\\ 

\begin{kor}\label{bflip} The functor $M\mapsto {\bf
  D}(\Theta_*{\mathcal V}_M)$ induces a bijection between

(a) the set of packets of standard supersingular ${\mathcal
  H}(G,I_0)_k$-modules, and

(b) the set of isomorphism classes of $B$-symmetric \'{e}tale
$(\varphi^{2d-1},\Gamma)$-modules ${\bf D}$ over $k_{\mathcal E}$.  
\end{kor}

{\sc Proof:} For $x\in{\mathbb F}_p^{\times}$ we compute$$\tau(x)\cdot \dot{u}\tau^{-1}(x)\dot{u}^{-1}={\rm
  diag}(x^2,E_{d-1},x^{-2},E_{d-1},1)=(\alpha^{\vee}_1-\alpha^{\vee}_0)(x)$$in
$\overline{T}$. Application of $\lambda$ gives $x^{s_f-s_e}=x^{k_1-k_0}$ and
hence $s_f-s_e\equiv k_1-k_0=k_{\beta(1)}-k_{\beta(2d-1)}$ modulo $(p-1)$. The
required symmetry in the $p$-adic digits of $n_e$, $n_f$ is due to the
corresponding symmetry of the function $\beta$. Thus, ${\bf
  D}(\Theta_*{\mathcal V}_M)$ is a $B$-symmetric \'{e}tale
$(\varphi^{2d-1},\Gamma)$-module. 

One checks that composing $M\mapsto{\bf
  D}(\Theta_*{\mathcal V}_M)$ with the map of Lemma \ref{combiflip} (ii) gives the map of Lemma \ref{combisusib}. As the maps in Lemmata \ref{combisusib} and \ref{combiflip} (ii) are bijective we obtain our result.\hfill$\Box$\\ 

\subsection{Affine root system $\tilde{D}_d$} 

\label{subsecd}

Assume $d\ge4$. Here $W_{{\rm aff}}$ is the
Coxeter group with generators $s_0, s_1,\ldots,s_d$ and relations
\begin{gather}(s_{d-2}s_d)^3=(s_{0}s_2)^3=(s_{i-1}s_i)^3=1\quad\mbox{ for }2\le i\le d-1\label{drel}\end{gather}and moreover
$(s_is_j)^2=1$ for all other pairs $i<j$, and $s_i^2=1$ for all $i$. In the extended affine Weyl group $\widehat{W}$
we find (cf. \cite{im}) an element $u$ of length $0$ with \begin{gather}u^2=1\quad\quad\quad\mbox{ and
}\quad\quad\quad us_0u=s_{1},\quad us_1u=s_{0},\quad us_{d-1}u=s_{d},\quad us_du=s_{d-1}\notag\\us_iu=s_{i}\quad\mbox{ for }2\le i\le d-2.\label{udfo}\end{gather}
(We have $\widehat{W}=W_{{\rm aff}}\rtimes W_{\Omega}$ with a four-element subgroup
$W_{\Omega}$, with $u\in W_{\Omega}$ of
order $2$.) 
 Consider the general orthogonal group $${\rm GO}_{2d}({\mathbb Q}_p)=\{A\in {\rm GL}_{2d}({\mathbb
  Q}_p)\,|\,{}^TAS_dA=\kappa(A)S_d\mbox{ for some }\kappa(A)\in{\mathbb
  Q}_p^{\times}\}.$$It contains the special orthogonal group $${\rm SO}_{2d}({\mathbb Q}_p)=\{A\in {\rm SL}_{2d}({\mathbb
  Q}_p)\,|\,{}^TAS_dA=S_d\}.$$

Let $G={\rm GSO}_{2d}({\mathbb Q}_p)$ be the subgroup of ${\rm
  GO}_{2d}({\mathbb Q}_p)$ generated by ${\rm
  SO}_{2d}({\mathbb Q}_p)$ and by all ${\rm diag}(xE_d,E_d)$ with $x\in{\mathbb
  Q}_p^{\times}$; it is of index $2$ in ${\rm
  GO}_{2d}({\mathbb Q}_p)$.\footnote{With the obvious definitions, ${\rm
  GO}_{2d}({\mathbb Q}_p)$ is the group of ${\mathbb Q}_p$-valued points of an algebraic group ${\rm
  GO}_{2d}$, whereas $G$ is the group of ${\mathbb Q}_p$-valued points of the connected component ${\rm GSO}_{2d}$ of ${\rm
  GO}_{2d}$; this ${\rm GSO}_{2d}$ has connected center.}  

Let $T$ be the maximal torus consisting of all diagonal matrices in
$G$. For $1\le i\le d$ let$$e_i:T\cap {\rm SL}_{2d}({\mathbb Q}_p)\longrightarrow {\mathbb
   Q}_p^{\times},\quad {\rm
  diag}(x_1,\ldots,x_{2d})\mapsto x_i.$$For $1\le i,j\le d$ and $\epsilon_1, \epsilon_2\in\{\pm 1\}$ we thus obtain characters (using additive notation as usual) $\epsilon_1e_i+\epsilon_2e_j:T\cap {\rm SL}_{2d}({\mathbb Q}_p)\longrightarrow {\mathbb
   Q}_p^{\times}$. We extend these latter ones to $T$ by setting$$\epsilon_1e_i+\epsilon_2e_j:T\longrightarrow {\mathbb
   Q}_p^{\times},\quad A={\rm
  diag}(x_1,\ldots,x_{2d})\mapsto
x_i^{\epsilon_1}x_j^{\epsilon_2}\kappa(A)^{\frac{-\epsilon_1-\epsilon_2}{2}}.$$Then $\Phi=\{\pm
e_i\pm e_j\,|\,i\ne j\}$ is the root system of $G$ with respect to $T$. It is
of type $D_d$. Choose the positive system $\Phi^+=\{
e_i\pm e_j\,|\,i< j\}$ with corresponding set of
simple roots $\Delta=\{\alpha_1=e_1-e_2,
\alpha_2=e_2-e_3,\ldots,\alpha_{d-1}=e_{d-1}-e_d, \alpha_{d}=e_{d-1}+e_d\}$. The negative of the highest root is $\alpha_0=-e_1-e_2$. For $0\le i\le d$ we have the following explicit formula for $\alpha_i^{\vee}=(\alpha_i)^{\vee}$:\begin{gather}\alpha_i^{\vee}(x)=\left\{\begin{array}{l@{\quad:\quad}l}{\rm diag}(x^{-1},x^{-1},E_{d-2},x,x,E_{d-2})&\quad
      i=0\\
{\rm
  diag}(E_{i-1},x,x^{-1},E_{d-i-1},E_{i-1},x^{-1},x,E_{d-i-1})&\quad 1\le
i\le d-1\\ 
{\rm diag}(E_{d-2},x,x,E_{d-2},x^{-1},x^{-1})&\quad
i=d\end{array}\right.\label{hsid}\end{gather}

For $\alpha\in \Phi$ let $N_{\alpha}^0$ be the subgroup of the
corresponding root subgroup $N_{\alpha}$ of $G$ all of whose elements belong
to ${\rm SL}_{2d}({\mathbb
  Z}_p)$. Let $I_0$ denote the pro-$p$-Iwahori subgroup generated by the $N_{\alpha}^0$ for
all $\alpha\in \Phi^+$, by the $(N_{\alpha}^0)^p$ for
all $\alpha \in \Phi^-=\Phi-\Phi^+$, and by the maximal pro-$p$-subgroup of $T_0$. Let $I$ denote the Iwahori
subgroup of $G$ containing $I_0$. Let $N_0$ be the subgroup
of $G$ generated by all $N_{\alpha}^0$ for $\alpha\in \Phi^+$. 

For $1\le i\le d-1$ define the block diagonal
matrix $$\dot{s}_i={\rm
  diag}(E_{i-1},{S}_1,E_{d-i-1},E_{i-1},{S}_1,E_{d-i-1})={\rm
  diag}(E_{i-1},{S}_1,E_{d-2},{S}_1,E_{d-i-1}).$$Put$$\dot{u}=\left(\begin{array}{cccccc}{}&{}&{}&p^{-1}&{}&{}\\{}&E_{d-2}&{}&{}&{}&{}\\{}&{}&{}&{}&{}&1\\p&{}&{}&{}&{}&{}\\{}&{}&{}&{}&E_{d-2}&{}\\{}&{}&1&{}&{}&{}\end{array}\right)$$and $\dot{s}_0=\dot{u}\dot{s}_1\dot{u}$ and
$\dot{s}_d=\dot{u}\dot{s}_{d-1}\dot{u}$. Then
$\dot{s}_{i}$ for $0\le i\le d$ belongs to ${\rm SO}_{2d}({\mathbb Q}_p)\subset G$ and normalizes $T$. Its image element $s_{i}=s_{\alpha_i}$ in $N(T)/ZT_0=\widehat{W}$ is the reflection corresponding to $\alpha_i$. The ${s}_0,{s}_1,\ldots,{s}_{d-1},{s}_{d}$ are
Coxeter generators of $W_{\rm aff}\subset \widehat{W}$ satisfying the
relations (\ref{drel}). Also $\dot{u}$ belongs to ${\rm SO}_{2d}({\mathbb Q}_p)\subset G$; it normalizes $T$, $I$ and $I_0$. Its image element $u$ in $N(T)/ZT_0=\widehat{W}$ satisfies the relations (\ref{udfo}). In $N(T)$ we consider\begin{align}\phi&=(p\cdot{\rm id})\dot{s}_{d-1}\dot{s}_{d-2}\cdots
\dot{s}_2\dot{s}_1\dot{s}_{d}\dot{s}_{d-2}\dot{s}_{d-3}\cdots
\dot{s}_3\dot{s}_2\dot{s}_0&\mbox{ if }d\mbox{ is even},\notag\\\phi&=(p^2\cdot{\rm id})\dot{s}_{d-1}\dot{s}_{d-2}\cdots
\dot{s}_2\dot{s}_1\dot{s}_{d}\dot{s}_{d-2}\dot{s}_{d-3}\cdots
\dot{s}_3\dot{s}_2\dot{s}_0&\mbox{ if }d\mbox{ is odd}.\notag\end{align}(The reason for our distinction according to the parity of $d$ is that the center $Z$ of $G$ is the subgroup of $G$ generated by $T_0\cap Z$ and $p\cdot{\rm id}$ if $d$ is even, resp. by $T_0\cap Z$ and $p^2\cdot{\rm id}$ if $d$ is odd.) We may rewrite this as $\phi=(p\cdot{\rm id})\dot{s}_{\beta(1)}\cdots
\dot{s}_{\beta(2d-2)}$ if $d$ is even, resp. $\phi=(p^2\cdot{\rm id})\dot{s}_{\beta(1)}\cdots
\dot{s}_{\beta(2d-2)}$ if $d$ is odd, where $\beta(i)=d-i$ for $1\le i\le d-1$, $\beta(d)=d$, $\beta(i)=2d-1-i$ for $d+1\le i\le 2d-3$ and $\beta(2d-2)=0$. We
put$$C^{(a(2d-2)+b)}=\phi^as_{\beta(1)}\cdots s_{\beta(b)}C$$for $a,b\in{\mathbb Z}_{\ge0}$ with $0\le b<2d-2$. Define the homomorphism$$\tau:{\mathbb
  Z}_p^{\times}\longrightarrow T_0,\quad\quad x\mapsto {\rm
  diag}(xE_{d-1},E_{d},x).$$

\begin{lem}\label{taudneu} We have $\phi^{d}\in T$ and $\phi^{d}N_0\phi^{-d}\subset N_0$. The sequence $C=C^{(0)}, C^{(1)},
C^{(2)},\ldots$ satisfies hypothesis (\ref{n0ineu}). In particular we may define $\alpha^{(j)}\in\Phi^+$
for all $j\ge0$. 

(b) For any $j\ge0$ we have $\alpha^{(j)}\circ\tau={\rm id}_{{\mathbb
  Z}_p^{\times}}$. 

(c) We have $\tau(a)\phi=\phi\tau(a)$ for all
$a\in{\mathbb Z}_p^{\times}$.
\end{lem}

{\sc Proof:} (a) A matrix computation shows $\phi^{d}={\rm
  diag}(p^{d+2}E_{d-1},p^{d-2}E_{d},p^{d+2})$ if $d$ is even, and $\phi^{d}={\rm
  diag}(p^{2d+2}E_{d-1},p^{2d-2}E_{d},p^{2d+2})$ if $d$ is odd. The group $N_{\alpha}$ for $\alpha\in\Phi^+$ is generated by $\epsilon_{i,j+d}\epsilon^{-1}_{j,i+d}$ if $\alpha=e_i+e_j$ with $1\le i<j\le d$, and by $\epsilon_{i,j}\epsilon^{-1}_{j+d,i+d}$ if $\alpha=e_i-e_j$ with $1\le i<j\le d$. Using this we
find$$\phi^{d}N_0\phi^{-d}=\prod_{\alpha\in\Phi^+}\phi^{d}(N_0\cap
N_{\alpha})\phi^{-d}=\prod_{\alpha\in\Phi^+}(N_0\cap
N_{\alpha})^{p^{m_{\alpha}}},$$$$m_{\alpha}=\left\{\begin{array}{l@{\quad:\quad}l}4&\quad
      \alpha=e_i+e_j\mbox{ with }1\le i<j<d\\
4&\quad\alpha=e_i-e_d\mbox{ with }1\le i<d\\
0&\quad\mbox{all other }\alpha\in\Phi^+\end{array}\right.$$In particular we find $\phi^{d}N_0\phi^{-d}\subset N_0$ and
$[N_0:\phi^{d}N_0\phi^{-d}]=p^{2d(d-1)}$. On the other hand, the image of $\phi$ in $\widehat{W}$ is a product of $2d-2$ Coxeter generators and of an element of length
$0$. Arguing as in the proof of Lemma \ref{tauc} we combine these facts to obtain $\ell(\phi)=2d-2$ and that $\phi$ is power multiplicative. We obtain that hypothesis
(\ref{n0ineu}) holds true, again by the same reasoning as in Lemma \ref{tauc}.

(b) As $\phi^{d}\in T$ we have
$\{\alpha^{(j)}\,|\,j\ge0\}=\{\alpha\in\Phi^+\,|\,m_{\alpha}\ne0\}$. This
implies (b). 


(c) Another matrix computation.\hfill$\Box$\\

As explained in subsection \ref{erklaer} we now obtain a functor from ${\rm Mod}^{\rm fin}({\mathcal
  H}(G,I_0))$ to the category of
$(\varphi^{2d-2},\Gamma)$-modules over ${\mathcal O}_{\mathcal E}$. Consider the
elements $$\dot{\omega}=\left(\begin{array}{cc}{}&E_{d}^*\\pE_{d}^*&{}\end{array}\right),\quad\quad\quad\dot{\rho}=\left(\begin{array}{cccc}{}&{}&{}&E_{d-1}^*\\p{}&{}&{}&{}\\{}&pE_{d-1}^*&{}&{}\\{}&{}&1&{}\end{array}\right)$$of
${\rm GO}_{2d}({\mathbb Q}_p)$. They normalize
$T$ and satisfy \begin{gather}\dot{\omega}\dot{u}=\dot{u}\dot{\omega},\quad\quad\dot{\omega}^2=p\cdot{\rm id} \notag\\\dot{\omega}\dot{s}_i\dot{\omega}^{-1}=\dot{s}_{d-i}\quad\mbox{ for }0\le i\le
d,\notag\\\dot{\rho}^2=p\cdot\dot{u},\notag\\\dot{\rho}\dot{s}_i\dot{\rho}^{-1}=\dot{s}_{d-i}\quad\mbox{ for }2\le i\le
d-2,\notag\\\dot{\rho}\dot{s}_{d-1}\dot{\rho}^{-1}=\dot{s}_{1},\quad\quad\dot{\rho}\dot{s}_{d}\dot{\rho}^{-1}=\dot{s}_0,\quad\quad\dot{\rho}\dot{s}_{0}\dot{\rho}^{-1}=\dot{s}_{d-1},\quad\quad\dot{\rho}\dot{s}_{1}\dot{\rho}^{-1}=\dot{s}_d.\notag\end{gather}

The element $\dot{\omega}$ belongs to $G$ if and only if $d$ is even. The
element $\dot{\rho}$ belongs to $G$ if and only if $d$ is odd.\\

 Let ${\mathcal H}(G,I_0)'_{{\rm
    aff},k}$ denote the $k$-sub algebra of ${\mathcal H}(G,I_0)_{k}$ generated
by ${\mathcal H}(G,I_0)_{{\rm
    aff},k}$ together with $T_{p\cdot {\rm id}}=T_{\dot{\omega}^{2}}$ and
$T^{-1}_{p\cdot {\rm id}}=T_{p^{-1}\cdot {\rm id}}$ if $d$ is even, resp. $T_{p^2\cdot {\rm id}}=T_{\dot{\rho}^{4}}$ and
$T^{-1}_{p^2\cdot {\rm id}}=T_{p^{-2}\cdot {\rm id}}$ if
$d$ is odd.

Suppose we are given a character $\lambda:\overline{T}\to
k^{\times}$, a
subset ${\mathcal J}\subset S_{\lambda}$ and some $b\in k^{\times}$. Define the numbers $0\le k_i=k_i(\lambda,{\mathcal
  J})\le p-1$ as in subsection \ref{basic}. The character $\chi_{\lambda,{\mathcal
  J}}$ of ${\mathcal H}(G,I_0)_{{\rm
    aff},k}$ extends uniquely to a character$$\chi_{\lambda,{\mathcal
  J},b}:{\mathcal H}(G,I_0)'_{{\rm
    aff},k}\longrightarrow k$$which sends $T_{p\cdot {\rm id}}$ to $b$ if $d$
is even, resp. which sends $T_{p^2\cdot {\rm id}}$ to $b$ if $d$ is odd (see the
proof of \cite{vigneras} Proposition 3). We
define the ${\mathcal H}(G,I_0)_k$-module$$M=M[\lambda,{\mathcal
  J},b]={\mathcal H}(G,I_0)_k\otimes_{{\mathcal H}(G,I_0)'_{{\rm
    aff},k}}k.e$$where $k.e$ denotes the one dimensional $k$-vector space on
the basis element $e$, endowed with the action of ${\mathcal H}(G,I_0)'_{{\rm
    aff},k}$ by the character $\chi_{\lambda,{\mathcal
  J},b}$. As a $k$-vector space, $M$ has dimension $4$. A $k$-basis is
$e_0,e_1,f_0,f_1$ where we write $$e_0=1\otimes e,\quad f_0=T_{\dot{u}}\otimes
e,\quad e_1=T_{\dot{\omega}}\otimes e,\quad f_1=T_{\dot{u}\dot{\omega}}\otimes
e\quad\quad\mbox{ if }d\mbox{ is even},$$$$e_0=1\otimes e,\quad f_0=T_{\dot{u}}\otimes
e,\quad e_1=T_{\dot{\rho}}\otimes e,\quad f_1=T_{\dot{u}\dot{\rho}}\otimes
e\quad\quad\mbox{ if }d\mbox{ is odd}.$$

{\bf Definition:} We call an ${\mathcal
  H}(G,I_0)_k$-module {\it standard supersingular} if it is isomorphic with $M[\lambda,{\mathcal
  J},b]$ for some $\lambda, {\mathcal
  J}, b$ such that $k_{\bullet}\notin\{(0,\ldots,0), (p-1,\ldots,p-1)\}$.\\

For $2\le j\le d-2$ let $\widetilde{j}=j$, and furthermore let
$\widetilde{d-1}=d$ and $\widetilde{d}=d-1$ and $\widetilde{1}=0$ and
$\widetilde{0}=1$. Letting $\widetilde{\beta}=\widetilde{(.)}\circ\beta$ we then
have\begin{gather}\dot{u}\phi\dot{u}^{-1}=(p^n\cdot{\rm
    id})\dot{s}_{\widetilde{\beta}(1)}\cdots\dot{s}_{\widetilde{\beta}(2d-2)}\notag\end{gather}with
$n=1$ if $d$ is even, but $n=2$ if $d$ is odd. If $d$ is odd we consider in addition the following two maps $\gamma$
and $\delta$ from $[1,2d-2]$ to $[0,d]$. We put $\gamma(1)=1$, $\gamma(d-1)=d$,
$\gamma(d)=0$ and $\gamma(2d-2)=d-1$. We put $\delta(1)=0$, $\delta(d-1)=d-1$,
$\delta(d)=1$ and $\delta(2d-2)=d$. We put $\gamma(i)=\delta(i)=\beta(2d-2-i)$
for all $i\in[2,\ldots,d-2]\cup[d+1,\ldots,2d-3]$. We then
have\begin{gather}\dot{\varrho}\phi\dot{\varrho}^{-1}=(p^2\cdot{\rm id})\dot{s}_{\gamma(1)}\cdots\dot{s}_{\gamma(2d-2)},\quad\quad\dot{\varrho}^{-1}\phi\dot{\varrho}=(p^2\cdot{\rm id})\dot{s}_{\delta(1)}\cdots\dot{s}_{\delta(2d-2)}.\notag\end{gather}Fix $M[\lambda,{\mathcal
  J},b]$ for some $\lambda, {\mathcal
  J}, b$. Put $$n_{e_0}=\sum_{i=0}^{2d-3}k_{\beta(i+1)}p^i,\quad\quad
n_{f_0}=\sum_{i=0}^{2d-3}k_{\widetilde{\beta}(i+1)}p^i\quad\quad\mbox{ for
  any parity of
}d,$$$$n_{e_1}=\sum_{i=0}^{2d-3}k_{\beta(2d-2-i)}p^i,\quad\quad
n_{f_1}=\sum_{i=0}^{2d-3}k_{\widetilde{\beta}(2d-2-i)}p^i\quad\quad\mbox{ if }d\mbox{ is even},$$$$n_{e_1}=\sum_{i=0}^{2d-3}k_{\gamma(i+1)}p^i,\quad\quad
n_{f_1}=\sum_{i=0}^{2d-3}k_{\delta(i+1)}p^i\quad\quad\mbox{ if
}d\mbox{ is odd}.$$Let $0\le s_{e_0}, s_{f_0},
s_{e_1}, s_{f_1}\le p-2$ be such that for all $x\in{\mathbb
  F}_p^{\times}$ we have$$\lambda(\tau(x))=x^{-s_{e_0}},\quad\quad\lambda(\dot{u}\tau(x)\dot{u}^{-1})=x^{-s_{f_0}}\quad\quad\mbox{ for
  any parity of
}d,$$$$\lambda(\dot{\omega}\tau(x)\dot{\omega}^{-1})=x^{-s_{e_1}},\quad\quad\lambda(\dot{\omega}\dot{u}\tau(x)\dot{u}^{-1}\dot{\omega}^{-1})=x^{-s_{f_1}}\quad\quad\mbox{ if }d\mbox{ is even},$$$$\lambda(\dot{\rho}\tau(x)\dot{\rho}^{-1})=x^{-s_{e_1}},\quad\quad\lambda(\dot{\rho}\dot{u}\tau(x)\dot{u}^{-1}\dot{\rho}^{-1})=x^{-s_{f_1}}\quad\quad\mbox{ if
}d\mbox{ is odd}.$$Put $\varrho=k_0!k_1!k_{d-1}!k_d!\prod_{i=2}^{d-2}(k_i!)^2=\prod_{i=0}^{2d-3}(k_{\beta(i+1)}!)=\prod_{i=0}^{2d-3}(k_{\widetilde{\beta}(i+1)}!)$.\\

\begin{lem}\label{combisusid} The assignment $M[\lambda,{\mathcal
  J},b]\mapsto (n_{e_0}, s_{e_0}, b\varrho^{-1})$ induces a bijection
between the set of isomorphism classes of standard supersingular ${\mathcal
  H}(G,I_0)_k$-modules and ${\mathfrak S}_D(2d-2)$. 
\end{lem}

{\sc Proof:} We have
$\alpha^{\vee}_{0}(x)\alpha^{\vee}_{1}(x)\alpha^{\vee}_{{d-1}}(x)\alpha^{\vee}_d(x)\prod_{i=2}^{d-2}(\alpha_{i}^{\vee})^2(x)=1$
for all $x\in{\mathbb F}_p^{\times}$ (as can be seen e.g. from formula
(\ref{hsid})). This implies\begin{gather}k_0+k_1+k_{d-1}+k_d+2\sum_{i=2}^{d-2}k_i\equiv n_{e_0}\equiv
  n_{f_0}\equiv n_{e_1}\equiv
  n_{f_1}\equiv0\quad\mbox{ mod }(p-1).\label{perreld}\end{gather}It follows
from \cite{vigneras} Proposition 3 that $M[\lambda,{\mathcal
  J},b]$ and $M[\lambda',{\mathcal
  J}',b']$ are isomorphic if and only if $b=b'$ and the pair $(\lambda, {\mathcal
  J})$ is conjugate with the pair $(\lambda', {\mathcal
  J}')$ by means of $\dot{u}^{n}\dot{\omega}^m$ for some $n,m\in\{0,1\}$ (if
$d$ is even), resp. by means of $\dot{u}^{n}\dot{\rho}^m$ for some $n,m\in\{0,1\}$ (if
$d$ is odd). Under the map $M[\lambda,{\mathcal
  J},b]\mapsto (n_{e_0}, s_{e_0}, b\varrho^{-1})$, conjugation by $\dot{u}$ corresponds to the permutation
$\iota_0$ of $\widetilde{\mathfrak S}_D(2d-2)$, while conjugation by
$\dot{\omega}$, resp. by $\dot{\rho}$, corresponds to the permutation
$\iota_1$ of $\widetilde{\mathfrak S}_D(2d-2)$. We may thus proceed as in the proof of Lemma
\ref{combisusic} to see that our mapping is well defined and bijective.\hfill$\Box$\\

\begin{lem}\label{globda} For an appropriate choice of the isomorphism
  $\Theta$ we have, in $H_0(\overline{\mathfrak X}_+,\Theta_*{\mathcal V}_M)$,\begin{align}t^{n_{e_j}}\varphi^{2d-2}e_j&=\varrho b^{-1}e_j,\notag\\t^{n_{f_j}}\varphi^{2d-2}
  f_j&=\varrho b^{-1}f_j,\notag\\\gamma(x)e_j&=x^{-s_{e_j}}e_j,\notag\\\gamma(x)f_j&=x^{-s_{f_j}}f_j\notag\end{align}for
$x\in{\mathbb F}_p^{\times}$ and $j=0,1$. The action of $\Gamma_0$ on
$H_0(\overline{\mathfrak X}_+,\Theta_*{\mathcal V}_M)$ is trivial on the subspace $M$.
\end{lem} 

{\sc Proof:} As in Lemma \ref{globca}.\hfill$\Box$\\

\begin{kor}\label{dphiga} The \'{e}tale $(\varphi^{2d-2},\Gamma)$-module over
  $k_{\mathcal E}$ associated with $H_0(\overline{\mathfrak
    X}_+,\Theta_*{\mathcal V}_M)$ admits a $k_{\mathcal E}$-basis $g_{e_0}$,
  $g_{f_0}$, $g_{e_1}$, $g_{f_1}$ such that for both $j=0$ and $j=1$ we have\begin{gather}\varphi^{2d-2}g_{e_j}=b\varrho^{-1}t^{n_{e_j}+1-p^{2d-2}}g_{e_j}\notag\\\varphi^{2d-2}g_{f_j}=b\varrho^{-1}t^{n_{f_j}+1-p^{2d-2}}g_{f_j}\notag\\\gamma(x)(g_{e_j})-x^{s_{e_j}}g_{e_j}\in t\cdot k_{\mathcal E}^+\cdot g_{e_j}\notag\\\gamma(x)(g_{f_j})-x^{s_{f_j}}g_{f_j}\in t\cdot k_{\mathcal E}^+\cdot g_{f_j}\notag\end{gather} 
\end{kor}

{\sc Proof:} This follows from Lemma \ref{globda} as explained in \cite{dfun}
Lemma 6.4.\hfill$\Box$\\

\begin{kor} The functor $M\mapsto {\bf
  D}(\Theta_*{\mathcal V}_M)$ induces a bijection between

(a) the set of isomorphism classes of standard supersingular ${\mathcal
  H}(G,I_0)_k$-modules, and

(b) the set of isomorphism classes of $D$-symmetric \'{e}tale
$(\varphi^{2d-2},\Gamma)$-modules over $k_{\mathcal E}$.  
\end{kor}

{\sc Proof:} We first verify that ${\bf
  D}(\Theta_*{\mathcal V}_M)$ for $M=M[\lambda,{\mathcal
  J},b]$ is indeed $D$-symmetric. For this let ${\bf D}_{11}=\langle g_{e_0}\rangle$, ${\bf
  D}_{12}=\langle g_{f_0}\rangle$, ${\bf D}_{21}=\langle g_{e_1}\rangle$,
${\bf D}_{22}=\langle g_{f_1}\rangle$. Then $k_i({\bf
  D}_{11})=k_{\beta(i+1)}$ and $k_i({\bf
  D}_{12})=k_{\widetilde{\beta}(i+1)}$; moreover $k_i({\bf
  D}_{21})=k_{\beta(2d-d-i)}$ and $k_i({\bf
  D}_{22})=k_{\widetilde{\beta}(2d-d-i)}$ if $d$ is even, but $k_i({\bf
  D}_{21})=k_{\gamma(i+1)}$ and $k_i({\bf
  D}_{22})=k_{\delta(i+1)}$ if $d$ is odd.  

For the condition on $s_{f_0}-s_{e_0}=s({\bf D}_{12})-s({\bf D}_{11})$ we
compute$$\tau(x)\cdot\dot{u}\tau^{-1}(x)\dot{u}^{-1}={\rm
  diag}(x,E_{d-2},x^{-1},x^{-1},E_{d-2},x)=(\sum_{i=1}^{d-1}\alpha^{\vee}_i)(x),$$hence
application of $\lambda$ gives $x^{s_{f_0}-s_{e_0}}=x^{\sum_{i=1}^{d-1}k_i}$ and
hence $s_{f_0}-s_{e_0}\equiv \sum_{i=1}^{d-1}k_i=\sum_{i=0}^{d-2}k_i({\bf D}_{11})$ modulo $(p-1)$. The
condition on $s_{f_1}-s_{e_1}=s({\bf D}_{22})-s({\bf D}_{21})$ in case $d$ is even is exactly verified like the one for $s_{f_0}-s_{e_0}$ because
$\dot{\omega}\tau(x)\dot{\omega}^{-1}\cdot\dot{\omega}\dot{u}\tau^{-1}(x)\dot{u}^{-1}\dot{\omega}^{-1}=\tau(x)\cdot\dot{u}\tau^{-1}(x)\dot{u}^{-1}$. In
case $d$ is odd the computation is $$\dot{\rho}\tau(x)\dot{\rho}^{-1}\cdot\dot{\rho}\dot{u}\tau^{-1}(x)\dot{u}^{-1}\dot{\rho}^{-1}={\rm
  diag}(x,E_{d-2},x,x^{-1},E_{d-2},x^{-1})=(\alpha_{d}^{\vee}+\sum_{i=1}^{d-2}\alpha^{\vee}_i)(x),$$hence
$s_{f_1}-s_{e_1}\equiv k_d+\sum_{i=1}^{d-2}k_i=s({\bf D}_{12})-s({\bf
  D}_{11})+k_d-k_{d-1}=\sum_{i=0}^{d-2}k_i({\bf D}_{21})$ modulo $(p-1)$. 

 To see the
condition on $s_{e_1}-s_{e_0}=s({\bf D}_{21})-s({\bf D}_{11})$ in case $d$ is even we
compute\begin{align}\tau(x)\cdot\dot{\omega}\tau^{-1}(x)\dot{\omega}^{-1}&={\rm
  diag}(1,xE_{d-2},1,1,x^{-1}E_{d-2},1)\\{}&=(\frac{d-2}{2}\alpha^{\vee}_{d-1}+\frac{d-2}{2}\alpha^{\vee}_{d}+\sum_{i=2}^{d-2}(i-1)\alpha^{\vee}_i)(x),\notag\end{align}hence
application of $\lambda$ gives $x^{s_{e_1}-s_{e_0}}=x^{\frac{d-2}{2}k_{d-1}+\frac{d-2}{2}k_{d}+\sum_{i=2}^{d-2}(i-1)k_i}$ and
hence $s_{e_1}-s_{e_0}\equiv
\frac{d-2}{2}k_{d-1}+\frac{d-2}{2}k_{d}+\sum_{i=2}^{d-2}(i-1)k_i=\frac{d-2}{2}(k_{d-1}({\bf D}_{11})+k_{0}({\bf D}_{11}))+\sum_{i=2}^{d-2}(i-1)k_{d-i-1}({\bf D}_{11})$ modulo
$(p-1)$. If however $d$ is odd we
compute\begin{align}\tau(x)\cdot\dot{\rho}\tau^{-1}(x)\dot{\rho}^{-1}&={\rm
  diag}(1,xE_{d-2},x^{-1},1,x^{-1}E_{d-2},x)\notag\\{}&=(\frac{d-1}{2}\alpha^{\vee}_{d-1}+\frac{d-3}{2}\alpha^{\vee}_{d}+\sum_{i=2}^{d-2}(i-1)\alpha^{\vee}_i)(x),\notag\end{align}hence
application of $\lambda$ gives $x^{s_{e_1}-s_{e_0}}=x^{\frac{d-1}{2}k_{d-1}+\frac{d-3}{2}k_{d}+\sum_{i=2}^{d-2}(i-1)k_i}$ and
hence $s_{e_1}-s_{e_0}\equiv
\frac{d-1}{2}k_{d-1}+\frac{d-3}{2}k_{d}+\sum_{i=2}^{d-2}(i-1)k_i=\frac{d-3}{2}k_{\frac{r}{2}}({\bf D}_{11})+\frac{d-1}{2}k_{0}({\bf D}_{11})+\sum_{i=2}^{d-2}(i-1)k_{d-i-1}({\bf D}_{11})$ modulo
$(p-1)$. We have shown that ${\bf
  D}(\Theta_*{\mathcal V}_M)$ is $D$-symmetric.

One checks that composing $M\mapsto{\bf
  D}(\Theta_*{\mathcal V}_M)$ with the map of Lemma \ref{dflip} gives the map of Lemma \ref{combisusid}. As the maps in Lemmata \ref{combisusid} and \ref{dflip} are bijective we obtain our result.\hfill$\Box$\\ 

{\bf Remark:} Consider the subgroup $G'={\rm SO}_{2d}({\mathbb Q}_p)$ of $G$. If we replace the above $\tau$ by $\tau:{\mathbb
  Z}_p^{\times}\longrightarrow T_0$, $x\mapsto {\rm
  diag}(xE_{d-1},x^{-1}E_{d},x)$, and if we replace the above $\phi$ by $\phi=\dot{s}_{d-1}\dot{s}_{d-2}\cdots
\dot{s}_2\dot{s}_1\dot{s}_{d}\dot{s}_{d-2}\dot{s}_{d-3}\cdots
\dot{s}_3\dot{s}_2\dot{s}_0$, then everything in fact happens inside
$G'$, and there is no dichotomy between $d$ even or odd. We then have $\alpha^{(j)}\circ\tau={\rm id}^2_{{\mathbb
  Z}_p^{\times}}$ for all $j\ge0$. Let ${\rm Mod}_0^{\rm fin}{\mathcal
  H}(G',G'\cap I_0)$ denote the category of
finite-${\mathfrak o}$-length ${\mathcal
  H}(G',G'\cap I_0)$-modules on which $\tau(-1)$ (i.e. $T_{\tau(-1)}=T_{\tau(-1)^{-1}}$) acts trivially. For $M\in{\rm Mod}_0^{\rm fin}{\mathcal
  H}(G',G'\cap I_0)$ we obtain an action of $\lfloor{\mathfrak N}_0,\varphi^{2d-2},\Gamma^2\rfloor$ on
$H_0(\overline{\mathfrak X}_+,\Theta_*{\mathcal V}_M)$, where
$\Gamma^2=\{\gamma^2\,|\,\gamma\in\Gamma\}\subset\Gamma$. Correspondingly, we obtain a functor from ${\rm Mod}^{\rm fin}_0{\mathcal
  H}(G',G'\cap I_0)$ to the
category of $(\varphi^{2d-2},\Gamma^2)$-modules over ${\mathcal O}_{\mathcal
  E}$.\\

{\bf Remark:} Instead of the element $\phi\in N(T)$ used above we might also
work with the element $\dot{s}_{d-1}\cdots
\dot{s}_2\dot{s}_1\dot{u}$ of length $d-1$ (or products of this with elements of
$p^{\mathbb Z}\cdot{\rm id}$), keeping the same $C^{(\bullet)}$. This results
in a functor from ${\mathcal
  H}(G,I_0)$-modules to $(\varphi^{d-1},\Gamma)$-modules. Up to a factor in $p^{\mathbb Z}\cdot{\rm id}$,
the square of $\dot{s}_{d-1}\cdots
\dot{s}_2\dot{s}_1\dot{u}$ is the element $\phi$ used above.\\

{\bf Remark:} For the affine root system of type $D_d$ there are three co
minuscule fundamental coweights (cf. \cite{bourb} chapter 8, par 7.3]). The corresponding $\phi$'s for the other two choices are
longer.

\subsection{Affine root system $\tilde{A}_d$}

Assume $d\ge1$ and consider $G={\rm GL}_{d+1}({\mathbb
  Q}_p)$. Let $$\dot{u}=\left(\begin{array}{cc}{}&E_d\\p&{}\end{array}\right).$$For
$1\le i\le d$ let $$\dot{s}_i={\rm diag}(E_{i-1},S_1,E_{d-i})$$and let
$\dot{s}_0=\dot{u}\dot{s}_1\dot{u}^{-1}$. Let $T$ be the maximal torus
consisting of diagonal matrices. Let $\Phi^+$ be such that
$N=\prod_{\alpha\in\Phi^+}N_{\alpha}$ is the subgroup of upper triangular
unipotent matrices. Let $I_0$ be the subgroup consisting of elements in ${\rm
  GL}_{d+1}({\mathbb Z}_p)$ which are upper triangular modulo $p$. We put $$\phi=(p\cdot{\rm id})\dot{s}_d\cdots\dot{s}_0=(p\cdot{\rm
  id})\dot{s}_{\beta(1)}\cdots\dot{s}_{\beta(d+1)}$$where
$\beta(i)=d+1-i$ for $1\le i\le d+1$. For $a,b\in{\mathbb Z}_{\ge0}$ with $0\le b<d+1$ we
put$$C^{(a(d+1)+b)}=\phi^as_d\cdots s_{d-b+1}C=\phi^as_{\beta(1)}\cdots
s_{\beta(b)}C.$$We define the homomorphism$$\tau:{\mathbb
  Z}_p^{\times}\longrightarrow T_0,\quad\quad x\mapsto{\rm
  diag}(E_d,x^{-1}).$$The sequence $C=C^{(0)}, C^{(1)},
C^{(2)},\ldots$ satisfies hypothesis
(\ref{n0ineu}). The corresponding $\alpha^{(j)}\in\Phi^+$
for $j\ge0$ satisfy $\alpha^{(j)}\circ\tau={\rm id}_{{\mathbb
  Z}_p^{\times}}$, and we have $\tau(a)\phi=\phi\tau(a)$ for all
$a\in{\mathbb Z}_p^{\times}$. We thus obtain a functor from ${\rm Mod}^{\rm fin}({\mathcal
  H}(G,I_0))$ to the category of
$(\varphi^{d+1},\Gamma)$-modules over ${\mathcal O}_{\mathcal E}$.

Let ${\mathcal H}(G,I_0)'_{{\rm
    aff},k}$ denote the $k$-sub algebra of ${\mathcal H}(G,I_0)_{k}$ generated
by ${\mathcal H}(G,I_0)_{{\rm
    aff},k}$ together with $T_{p\cdot {\rm id}}=T_{\dot{u}^{d+1}}$ and
$T^{-1}_{p\cdot {\rm id}}=T_{p^{-1}\cdot {\rm id}}$.

Suppose we are given a character $\lambda:\overline{T}\to
k^{\times}$, a
subset ${\mathcal J}\subset S_{\lambda}$ and some $b\in k^{\times}$. Define the numbers $0\le k_i=k_i(\lambda,{\mathcal
  J})\le p-1$ as in subsection \ref{basic}. The character $\chi_{\lambda,{\mathcal
  J}}$ of ${\mathcal H}(G,I_0)_{{\rm
    aff},k}$ extends uniquely to a character$$\chi_{\lambda,{\mathcal
  J},b}:{\mathcal H}(G,I_0)'_{{\rm
    aff},k}\longrightarrow k$$which sends $T_{p\cdot {\rm id}}$ to $b$ (see the
proof of \cite{vigneras} Proposition 3). Define the ${\mathcal H}(G,I_0)_k$-module$$M=M[\lambda,{\mathcal
  J},b]={\mathcal H}(G,I_0)_k\otimes_{{\mathcal H}(G,I_0)'_{{\rm
    aff},k}}k.e$$where $k.e$ denotes the one dimensional $k$-vector space on
the basis element $e$, endowed with the action of ${\mathcal H}(G,I_0)'_{{\rm
    aff},k}$ by the character $\chi_{\lambda,{\mathcal
  J},b}$. As a $k$-vector space, $M$ has dimension $d+1$, a $k$-basis is
$\{e_i\}_{0\le i\le d}$ where we write $e_i=T_{\dot{u}^{-i}}\otimes e$.\\

{\bf Definition:} We call an ${\mathcal
  H}(G,I_0)_k$-module {\it standard supersingular} if it is isomorphic with $M[\lambda,{\mathcal
  J},b]$ for some $\lambda, {\mathcal
  J}, b$ such that $k_{\bullet}\notin\{(0,\ldots,0), (p-1,\ldots,p-1)\}$.\\

For $0\le j\le d$ put $n_{e_j}=\sum_{i=0}^dk_{j-i}p^i$ (reading $j-i$ as its representative modulo $(d+1)$ in $[0,d]$) and let $s_{e_j}$ be
such that $\lambda(\dot{u}^{-j}\tau(x)\dot{u}^j)=x^{-s_{e_j}}$. Put
$\varrho=\lambda(-{\rm id})\prod_{i=0}^d(k_i!)$. 

\begin{satz} The \'{e}tale $(\varphi^{d+1},\Gamma)$-module ${\bf
  D}(\Theta_*{\mathcal V}_M)$ over $k_{\mathcal E}$ associated with
$H_0(\overline{\mathfrak X}_+,\Theta_*{\mathcal V}_M)$ admits a $k_{\mathcal
  E}$-basis $\{g_{e_j}\}_{0\le j\le d}$ such that for all $j$ we have\begin{gather}\varphi^{d+1}g_{e_j}=b\varrho^{-1}t^{n_{e_j}+1-p^{d+1}}g_{e_j},\notag\\\gamma(x)g_{e_j}-x^{s_{e_j}}g_{e_j}\in t\cdot k_{\mathcal E}^+\cdot g_{e_j}.\notag\end{gather}The functor $M\mapsto {\bf
  D}(\Theta_*{\mathcal V}_M)$ induces a bijection between

(a) the set of isomorphism classes of standard supersingular ${\mathcal
  H}(G,I_0)_k$-modules, and

(b) the set of isomorphism classes of $A$-symmetric \'{e}tale
$(\varphi^{d+1},\Gamma)$-modules over $k_{\mathcal E}$.  
\end{satz}

{\sc Proof:} For the formulae describing ${\bf
  D}(\Theta_*{\mathcal V}_M)$ one may proceed exactly as in the proof of
Corollary \ref{cphigaa}. (The only tiny additional point to be observed is
that the $\dot{s}_i$ (in keeping with our choice in \cite{dfun}) do not ly in the images of the $\iota_{\alpha_i}$;
this is accounted for by the sign factor $\lambda(-{\rm id})$ in the
definition of $\varrho$.) Alternatively, as our $\phi$ is the $(d+1)$-st power
of the $\phi$ considered in section 8 of \cite{dfun}, the computations of
loc. cit. may be carried over. 

To see that ${\bf
  D}(\Theta_*{\mathcal V}_M)$ is $A$-symmetric put ${\bf D}_j=\langle
g_{e_j}\rangle$ for $0\le j\le d$ and compare the above formulae with those
defining $A$-symmetry; e.g. we find $s_{e_0}-s_{e_j}\equiv\sum_{i=1}^jk_i$
modulo $(p-1)$. The bijectivity statement is then verified as before.\hfill$\Box$\\

{\bf Remark:} Application of the functor of Lemma \ref{gamsemilin} to
any one of the direct summands ${\bf D}_j$ of the $A$-symmetric \'{e}tale
$(\varphi^{d+1},\Gamma)$-module ${\bf
  D}(\Theta_*{\mathcal V}_M)$ yields an \'{e}tale $(\varphi,\Gamma)$-module
isomorphic with the one assigned to $M$ in \cite{dfun}.\\

{\bf Remark:} Consider the subgroup $G'={\rm SL}_{d+1}({\mathbb Q}_p)$ of $G$. If we replace the above $\tau$ by $\tau:{\mathbb
  Z}_p^{\times}\longrightarrow T_0$, $x\mapsto {\rm
  diag}(xE_{d},x^{-d})$ and if we replace the above $\phi$ by $\phi=\dot{s}_{d}\dot{s}_{d-1}\cdots
\dot{s}_{1}\dot{s}_0$ then everything in fact happens inside
$G'$. We then have $\alpha^{(j)}\circ\tau={\rm id}^{d+1}_{{\mathbb
  Z}_p^{\times}}$ for all $j\ge0$. Let ${\rm Mod}_0^{\rm fin}{\mathcal
  H}(G',G'\cap I_0)$ denote the category of
finite-${\mathfrak o}$-length ${\mathcal
  H}(G',G'\cap I_0)$-modules on which the $xE_{d+1}$ (i.e. the $T_{x^{-1}E_{d+1}}$) for all $x\in{\mathbb Z}_p^{\times}$ with
$x^{d+1}=1$ act trivially. (Notice that $\tau(x)=xE_{d+1}$ for such $x$.) For $M\in{\rm Mod}_0^{\rm fin}{\mathcal
  H}(G',G'\cap I_0)$ we obtain an action of $\lfloor{\mathfrak N}_0,\varphi^{d+1},\Gamma^{d+1}\rfloor$ on
$H_0(\overline{\mathfrak X}_+,\Theta_*{\mathcal V}_M)$, where
$\Gamma^{d+1}=\{\gamma^{d+1}\,|\,\gamma\in\Gamma\}\subset\Gamma$. Correspondingly, we obtain a functor from ${\rm Mod}^{\rm fin}_0{\mathcal
  H}(G',G'\cap I_0)$ to the category of $(\varphi^{d+1},\Gamma^{d+1})$-modules over ${\mathcal O}_{\mathcal E}$.  \\

{\bf Remark:} As all the fundamental coweights $\tau$ of $T$ are minuscule, each of them
admits $\phi$'s for which the pair $(\phi,\tau)$ satisfies the properties
asked for in Lemma \ref{concept}. For example, let $1\le g\le d$. For $\phi=\dot{s}_g\cdot
\dot{s}_{g+1}\cdots \dot{s}_d\cdot \dot{u}$ as well as for $\phi=\dot{s}_g\cdot \dot{s}_{g-1}\cdots
\dot{s}_1\cdot \dot{u}^{-1}$ there is a unique minimal gallery from $C$ to $\phi(C)$ which
admits a $\phi$-periodic continuation to a gallery (\ref{cgall}),  giving rise
to a functor from ${\rm Mod}^{\rm fin}({\mathcal
  H}(G,I_0))$ to the category of
$(\varphi^{r},\Gamma)$-modules over ${\mathcal O}_{\mathcal E}$, where $r=\ell(\phi)$.


\section{Exceptional groups of type $\tilde{E}_6$ and $\tilde{E}_7$}

Let $G$ be the group of ${\mathbb Q}_p$-rational points of a ${\mathbb
  Q}_p$-split connected reductive group over ${\mathbb Q}_p$ with connected center $Z$. Fix a maximal
${\mathbb Q}_p$-split torus $T$ and define $\Phi$, $N(T)$, $W$, $\widehat{W}$ and $W_{{\rm aff}}$ as before.

\subsection{Affine root system $\tilde{E}_6$} 

Assume that the root system $\Phi$ is of type $E_6$. Following \cite{bourb} (for the indexing) we then have generators
$s_1,\ldots,s_6$ of $W$ and $s_0,\ldots,s_6$ of $W_{{\rm
    aff}}$ such that
$$(s_{1}s_3)^3=(s_{3}s_4)^3=(s_{4}s_5)^3=(s_{5}s_6)^3=(s_{2}s_4)^3=(s_{0}s_2)^3=1$$and moreover
$(s_is_j)^2=1$ for all other pairs $i<j$, and $s_i^2=1$ for all $i$. In the extended affine Weyl group $\widehat{W}$
we find (cf. \cite{im}) an element $u$ of length $0$ with \begin{gather}u^3=1\quad\quad\quad\mbox{ and
}\quad\quad\quad us_4u^{-1}=s_{4},\notag\\us_3u^{-1}=s_{5},\quad us_{5}u^{-1}=s_{2},\quad us_2u^{-1}=s_{3},\label{wufo}\\us_1u^{-1}=s_{6},\quad us_{6}u^{-1}=s_{0},\quad us_0u^{-1}=s_{1}.\notag\end{gather}
(We have $\widehat{W}=W_{{\rm aff}}\rtimes W_{\Omega}$ with the three-element subgroup $W_{\Omega}=\{1,u,u^2\}$.) Let $e_1,\ldots,e_8$ denote the standard basis of ${\mathbb R}^8$. We use the
standard inner product $\langle.,.\rangle$ on ${\mathbb R}^8$ to view both the
root system $\Phi$ as well as its dual $\Phi^{\vee}$ as living inside
${\mathbb R}^8$. We choose a positive system $\Phi^+$ in $\Phi$ such that, as in \cite{bourb}, the simple
roots are $\alpha_1=\alpha_1^{\vee}=\frac{1}{2}(e_1+e_8-e_2-e_3-e_4-e_5-e_6-e_7)$,
$\alpha_2=\alpha_2^{\vee}=e_2+e_1$, $\alpha_3=\alpha_3^{\vee}=e_2-e_1$, $\alpha_4=\alpha_4^{\vee}=e_3-e_2$,
$\alpha_5=\alpha_5^{\vee}=e_4-e_3$, $\alpha_6=\alpha_6^{\vee}=e_5-e_4$ while
the negative of the highest root is
$\alpha_0=\alpha_0^{\vee}=\frac{1}{2}(e_6+e_7-e_1-e_2-e_3-e_4-e_5-e_8)$. The set of positive roots is $$\Phi^+=\{e_j\pm e_i\,|\,1\le i<j\le 5\}\cup\{\frac{1}{2}(-e_6-e_7+e_8+\sum_{i=1}^5(-1)^{\nu_i}e_i)\,|\,\sum_{i=1}^5\nu_i\mbox{ even}\}.$$
We lift $u$ and the $s_i$ to elements $\dot{u}$ and $\dot{s}_i$ in $N(T)$. We then put \begin{gather}\phi=\dot{s}_2\dot{s}_4\dot{s}_3\dot{s}_1\dot{u}^{-1}\in N(T).\label{e6dua}\end{gather}We define $\nabla$ as in section \ref{powbas}.

\begin{pro}\label{e6comi}  There is a $\tau\in \nabla$ such that the pair $(\phi,\tau)$ satisfies the hypotheses of Lemma \ref{concept}. More precisely, $\phi$ is power multiplicative, and for the minuscule fundamental (co)weight $\tau=\omega_1=\frac{2}{3}(e_8-e_7-e_6)\in \nabla$ we have $\phi^{12}=\tau^3$ in $W_{{\rm aff}}$.
\end{pro}

{\sc Proof:} (Here $\tau^3$ designates the three fold iterate of translation by $\tau=\omega_1$, i.e. translation by the element $3\omega_1$ of the root lattice.) We define the subset $$\Phi(\omega_1)=\{\frac{1}{2}(-e_6-e_7+e_8+\sum_{i=1}^5(-1)^{\nu_i}e_i)\,|\,\sum_{i=1}^5\nu_i\mbox{
  even}\}$$ of $\Phi^+$. We compute $\langle\beta,3\omega_1\rangle=3$ for each $\beta\in \Phi(\omega_1)$, but $\langle\beta,3\omega_1\rangle=0$ for each $\beta\in \Phi^+-\Phi(\omega_1)$. This means that for each $\beta\in \Phi(\omega_1)$ the translation by $3\omega_1$ crosses $3$ walls parallel to the hyperplane corresponding to $\beta$, and that for each $\beta\in \Phi^+-\Phi(\omega_1)$ it crosses no wall parallel to the hyperplane corresponding to $\beta$. As $\Phi(\omega_1)$ contains 16 elements it follows that $\ell(\tau^3)=\ell(3\omega_1)=3\cdot 16=48$. On the other hand, in the appendix we explain how one can prove $\phi^{12}=3\omega_1$ (as elements in $W_{{\rm aff}}$). Combining these two facts we deduce $\ell(\phi^{12})=48$ and hence that $\phi$ (which a priori has length $\le 4$) is power multiplicative.\hfill$\Box$\\

{\bf Remark:} Instead of verifying $\phi^{12}=3\omega_1$ by means of a direct (computer) computation, one might try to prove that the affine transformation $\phi^{12}$ on $A$ crosses exactly the same walls as does $3\omega_1$ (which are described in the above proof).\\

As explained in subsection \ref{erklaer} we now obtain a functor from ${\rm Mod}^{\rm fin}({\mathcal
  H}(G,I_0))$ to the category of
$(\varphi^{4},\Gamma)$-modules over ${\mathcal O}_{\mathcal E}$. 

Similarly, we may replace $\phi$ by its third power $\phi^{3}$ which (in
contrast to $\phi$) is an element of $W_{\rm aff}$ (modulo
$T_0$). It yields a functor from ${\rm Mod}^{\rm fin}({\mathcal
  H}(G,I_0))$ to the category of
$(\varphi^{12},\Gamma)$-modules over ${\mathcal O}_{\mathcal E}$. As in our
treatment of the cases $C$, $B$, $D$ and $A$, this functor identifies the set of
standard supersingular ${\mathcal
  H}(G,I_0)_k$-modules bijectively with a set of certain $E$-symmetric
\'{e}tale $(\varphi^{12},\Gamma)$-modules over ${k}_{\mathcal E}$ of dimension
$3$. We leave the details to the reader.\\

{\bf Remarks:} (a) Dual to the above choice of $\phi$ is the choice\begin{gather}\phi=\dot{s}_2\dot{s}_4\dot{s}_5\dot{s}_6\dot{u}\in N(T).\label{e6du}\end{gather}For this choice, Proposition \ref{e6comi} holds true verbatim the same way, but now with the minuscule fundamental (co)weight $\tau=\omega_6=\frac{1}{3}(3e_5+e_8-e_7-e_6)$ with its corresponding subset (again containing $16$ elements)$$\Phi(\omega_6)=\{\frac{1}{2}(e_5-e_6-e_7+e_8+\sum_{i=1}^4(-1)^{\nu_i}e_i)\,|\,\sum_{i=1}^4\nu_i\mbox{
  even}\}\cup\{e_5\pm e_i\,|\,1\le i< 5\}$$of $\Phi^+$. Again see the appendix.

(b) For $\phi$ given by either (\ref{e6dua}) or (\ref{e6du}), consider the corresponding reduced expression of $\phi^3$ in $W_{{\rm aff}}$ (obtained by three fold concatenation and then conjugation of the $s_i$'s appearing with powers of $u$ so that no factor $u$ or $u^2$ remains). The number of occurencess of the $s_i$ are precisely the coefficients of the $\alpha^{\vee}_i$ in$$\alpha_0^{\vee}+\alpha^{\vee}_1+\alpha_6^{\vee}+2\alpha^{\vee}_2+2\alpha^{\vee}_3+2\alpha^{\vee}_5+3\alpha^{\vee}_4=0.$$

\subsection{Affine root system $\tilde{E}_7$} 

Assume that the root system $\Phi$ is of type $E_7$. Following \cite{bourb} we then have generators
$s_1,\ldots,s_7$ of $W$ and $s_0,\ldots,s_7$ of $W_{{\rm
    aff}}$ such that
$$(s_{0}s_1)^3=(s_{1}s_3)^3=(s_{3}s_4)^3=(s_{4}s_5)^3=(s_{5}s_6)^3=(s_{6}s_7)^3=(s_{2}s_4)^3=1$$and moreover
$(s_is_j)^2=1$ for all other pairs $i<j$, and $s_i^2=1$ for all $i$. In the extended affine Weyl group $\widehat{W}$
we find (cf. \cite{im}) an element $u$ of length $0$ with \begin{gather}u^2=1\quad\quad\quad\mbox{ and
}\quad\quad\quad us_4u=s_{4},\quad us_2u=s_{2},\notag\\ us_3u=s_{5},\quad us_{6}u=s_{1},\quad us_7u=s_{0},\notag\\us_5u=s_{3},\quad us_{1}u=s_{6},\quad us_0u=s_{7}.\notag\end{gather}
(We have $\widehat{W}=W_{{\rm aff}}\rtimes W_{\Omega}$ with the two-element subgroup $W_{\Omega}=\{1,u\}$.) Let $e_1,\ldots,e_8$ denote the standard basis of ${\mathbb R}^8$. We use the standard inner product $\langle.,.\rangle$ on ${\mathbb R}^8$ to view both the root system $\Phi$ as well as its dual $\Phi^{\vee}$ as living inside ${\mathbb R}^8$. We choose a positive system $\Phi^+$ in $\Phi$ such that, as in \cite{bourb}, the simple
roots are $\alpha_1=\alpha_1^{\vee}=\frac{1}{2}(e_1+e_8-e_2-e_3-e_4-e_5-e_6-e_7)$,
$\alpha_2=\alpha_2^{\vee}=e_2+e_1$, $\alpha_3=\alpha_3^{\vee}=e_2-e_1$, $\alpha_4=\alpha_4^{\vee}=e_3-e_2$,
$\alpha_5=\alpha_5^{\vee}=e_4-e_3$, $\alpha_6=\alpha_6^{\vee}=e_5-e_4$, $\alpha_7=\alpha_7^{\vee}=e_6-e_5$ while
the negative of the highest root is
$\alpha_0=\alpha_0^{\vee}=e_7-e_8$. The set of positive roots is$$\Phi^+=\{e_j\pm e_i\,|\,1\le i<j\le 6\}\cup\{e_8-e_7\}\cup\{\frac{1}{2}(e_8-e_7+\sum_{i=1}^6(-1)^{\nu_i}e_i)\,|\,\sum_{i=1}^6\nu_i\mbox{ odd}\}.$$We lift $u$ and the $s_i$ to elements $\dot{u}$ and $\dot{s}_i$ in $N(T)$. We then put $$\phi=\dot{s}_1\dot{s}_3\dot{s}_4\dot{s}_2\dot{s}_5\dot{s}_4\dot{s}_3\dot{s}_1\dot{s}_0\dot{u}\in N(T).$$We define $\nabla$ as in section \ref{powbas}.

\begin{pro}\label{e7comi}  There is a $\tau\in \nabla$ such that the pair $(\phi,\tau)$ satisfies the hypotheses of Lemma \ref{concept}. More precisely, $\phi$ is power multiplicative, and for the minuscule fundamental (co)weight $\tau=\omega_7=e_6+\frac{1}{2}(e_8-e_7)\in \nabla$ we have $\phi^{6}=\tau^2$ in $W_{{\rm aff}}$.
\end{pro}

{\sc Proof:} Exactly the same as for Propostion \ref{e6comi}. The
corresponding subset in $\Phi^+$ is$$\{\frac{1}{2}(e_6+e_8-e_7+\sum_{i=1}^5(-1)^{\nu_i}e_i)\,|\,\sum_{i=1}^5\nu_i\mbox{
  odd}\}\cup\{e_8-e_7\}\cup\{e_6\pm e_i\,|\,1\le i< 6\}.$$It contains 27 elements,
thus $\ell(2\omega_7)=54$. For a computer proof of $\phi^{6}=2\omega_7$ see the appendix.\hfill$\Box$\\ 

As explained in subsection \ref{erklaer} we now obtain a functor from ${\rm Mod}^{\rm fin}({\mathcal
  H}(G,I_0))$ to the category of
$(\varphi^{9},\Gamma)$-modules over ${\mathcal O}_{\mathcal E}$. 

Similarly, we may replace $\phi$ by its square $\phi^{2}$ which (in
contrast to $\phi$) is an element of $W_{\rm aff}$ (modulo
$T_0$). It yields a functor from ${\rm Mod}^{\rm fin}({\mathcal
  H}(G,I_0))$ to the category of
$(\varphi^{18},\Gamma)$-modules over ${\mathcal O}_{\mathcal E}$. Again this functor identifies the set of
standard supersingular ${\mathcal
  H}(G,I_0)_k$-modules bijectively with a set of certain $E$-symmetric
\'{e}tale $(\varphi^{18},\Gamma)$-modules over ${k}_{\mathcal E}$ of dimension
$2$. We leave the details to the reader.\\

{\bf Remark:} The number of occurencess of the $s_i$ in $\phi^2\in W_{{\rm aff}}$ (cf. the case $\tilde{E}_6$) are the coefficients of the $\alpha^{\vee}_i$ in$$\alpha_0^{\vee}+\alpha^{\vee}_7+2(\alpha_1^{\vee}+\alpha_2^{\vee}+\alpha_6^{\vee})+3(\alpha^{\vee}_3+\alpha^{\vee}_5)+4\alpha^{\vee}_4=0.$$

\section{Appendix}

{\bf Verification of the statement $\phi^{12}=3\omega_1$ in the proof of Proposition \ref{e6comi}.}\\In the computer algebra system {\it sage}, the input

{\it R=RootSystem(["E",6,1]).weight\_lattice()

Lambda=R.fundamental\_weights()

omega1=Lambda[1]-Lambda[0]

R.reduced\_word\_of\_translation(3*omega1)}\\prompts the output\begin{gather}[0, 2, 4, 3, 5, 4, 2, 0, 6, 5, 4, 2, 3, 1, 4, 3, 5, 4, 2, 0, 6, 5, 4, 2,\notag\\3, 1, 4, 3, 5, 4, 2, 0, 6, 5, 4, 2, 
3, 1, 4, 3, 5, 4, 2, 6, 5, 4, 3, 1]\label{sagehelp}.\end{gather}By definition
of the function {\it reduced\_word\_of\_translation} this means $s^*_{i_1}\cdots s^*_{i_{48}}=3\omega_1$, with the string $[i_1,\ldots,i_{48}]$ as given by
(\ref{sagehelp}). Here $s_i^*=s_i$ for $1\le i\le 6$, but $s_0^*$ denotes the affine reflection
in the outer face of Bourbaki's fundamental alcove $A$. Since we deviate from
these conventions in that our $s_0$ is the affine reflection in the outer face
of the {\it negative} $C=-A$ of $A$, we must modify the above string
(\ref{sagehelp}) as follows. First, writing
$s_i^{**}=s_0^*s_i^*s_0^*$ for $0\le i\le 6$, conjugating the factors in the
previous word by $s_0^*$ and commuting some of its factors where allowed, the above says $s^{**}_{j_1}\cdots
s^{**}_{j_{48}}=3\omega_1$ where the string $[j_1,\ldots,j_{48}]$ is given
by\begin{gather}[2, 4, 5, 6, 3, 4, 2, 0, 5, 4, 3, 1, 2, 4, 5, 6, 3, 4, 2, 0,
  5, 4, 3, 1,\notag\\2, 4, 5, 6, 3, 4, 2, 0, 5, 4, 3,
  1, 2, 4, 5, 6, 3, 4, 2, 0, 5, 4, 3,
  1]\label{stsagehelp}.\end{gather}The $s_i^{**}$ are precisely the
reflections in the codimension $1$ faces of $s_0^*A$. But $s_0^*A$ is a
translate of $C$, and under this translation, the reflection $s_0^{**}=s_0^*$
corresponds to $s_0$, whereas for $1\le i\le 6$ the reflection $s_i^{**}$ corresponds to $w_0s_iw_0$, where
$w_0$ is the longest element of $W$. We have $w_0s_iw_0=s_i$
for $i\in\{0,2,4\}$, but $w_0s_3w_0=s_5$ and $w_0s_1w_0=s_6$. Thus, we obtain $s_{k_1}\cdots
s_{k_{48}}=3\omega_1$ where the string $[k_1,\ldots,k_{48}]$ is obtained
from the string (\ref{stsagehelp}) by keeping its entry values $0, 2$
and $4$, while exchanging the entry values $3$ with $5$ and $1$ with $6$. Using formulae
(\ref{wufo}) one checks that $s_{k_1}\cdots s_{k_{48}}=\phi^{12}$.\hfill$\Box$\\

{\bf Verification of the statement $\phi^{12}=3\omega_6$ for $\phi$ given by (\ref{e6du}).}\\The argument is the same as in Proposition \ref{e6comi}. The string returned by {\it sage} to 

{\it R=RootSystem(["E",6,1]).weight\_lattice(),} {\it
  Lambda=R.fundamental\_weights(),} 

{\it omega6=Lambda[6]-Lambda[0],} {\it
  R.reduced\_word\_of\_translation(3*omega6)}\\reads\begin{gather}[0, 2, 4, 3, 1, 5, 4, 2, 0, 3, 4, 2, 5, 4, 3, 1, 6, 5, 4, 2, 0, 3, 4, 2,\notag\\5, 4, 3, 1, 6, 5, 4, 2, 0, 3, 4, 2, 
5, 4, 3, 1, 6, 5, 4, 2, 3, 4, 5, 6].\notag\end{gather}

{\bf Verification of the statement $\phi^{6}=2\omega_7$ in the proof of Proposition \ref{e7comi}.}\\The argument is the same as in Proposition \ref{e6comi}. The string returned by {\it sage} to 

{\it
  R=RootSystem(["E",7,1]).weight\_lattice(),} {\it Lambda=R.fundamental\_weights(),} 

{\it omega7=Lambda[7]-Lambda[0],} {\it R.reduced\_word\_of\_translation(2*omega7)}\\reads\begin{gather}[0, 1, 3, 4,
  2, 5, 4, 3, 1, 0, 6, 5, 4, 2, 3, 1, 4, 3, 5, 4, 2, 6, 5, 4, 3, 1, 0,\\7, 6,
  5, 4, 2, 3, 1, 4, 3, 5, 4, 2, 6, 5, 4, 3, 1, 7, 6, 5, 4, 2, 3, 4, 5, 6,
  7].\notag\end{gather}

\begin{flushleft} \textsc{Humboldt-Universit\"at zu Berlin\\Institut f\"ur Mathematik\\Rudower Chaussee 25\\12489 Berlin, Germany}\\ \textit{E-mail address}:
gkloenne@math.hu-berlin.de \end{flushleft} \end{document}